\documentclass[oneside,12pt]{amsart}

\usepackage{stddefs}
\usepackage{defs-more}
\usepackage[matrix,arrow,curve]{xy}
\usepackage{hyperref}

\numberwithin{equation}{section}

\usepackage[top=1in,bottom=1in,left=1in,right=1in]{geometry}

\copyrightinfo{2004}{Andrew Mauer-Oats}

\begin{document}


\title[Algebraic Goodwillie Calculus]{Algebraic Goodwillie Calculus
  \\ and a Cotriple Model for the Remainder}

\author{Andrew Mauer-Oats}
\address{Department of Mathematics \\ Purdue University \\ 
West Lafayette, Indiana 
47907}
\email{amauer@math.purdue.edu}

\date{\today}



\begin{abstract}
Goodwillie has defined a tower of approximations for a functor from
spaces to spaces that is analogous to the Taylor series of a
function. His \nth{} order approximation $P_n F$ at a space $X$
depends on the values of $F$ on coproducts of large suspensions of the
space: $F(\vee \Sigma^M X)$.

We define an ``algebraic'' version of the Goodwillie tower, $\Pnalg
F(X)$, that depends only on the behavior of $F$ on coproducts of $X$.
When $F$ is a functor to connected spaces or grouplike $H$-spaces, the
functor $\Pnalg F$ is the base of a fibration 
$$ 
\realization{\Perp^{*+1} F}
\rightarrow 
F
\rightarrow 
\Pnalg F,
$$
whose fiber is the simplicial space associated to a cotriple $\Perp$
built from the $(n+1)^{\text{st}}$ cross effect of the functor $F$.
In a range in which $F$ commutes with realizations (for instance,
when $F$ is the identity functor of spaces),
the algebraic Goodwillie tower agrees with the ordinary (topological)
Goodwillie tower, so this theory gives a way of studying the Goodwillie
approximation to a functor $F$ in many interesting cases. 
\end{abstract}

\maketitle


\section{Introduction}


A function on the real line can be approximated by its ordinary Taylor
series at a point, with the ``\nth{}-order approximation'' being the
Taylor polynomial through the $x^n$ term.  Goodwillie
\cite{Cal1,Cal2,Cal3} has defined an analogous approximation for
functors of topological spaces.

Johnson and McCarthy
\cite{Johnson-McCarthy:taylor-towers-for-functors-of-additive-categories}
have explored alternative models for Goodwillie's approximation,
working in the simpler setting of chain complexes.
Specifically,
for a functor $C$ that is the prolongation to chain complexes of a
functor to an abelian category, they
have defined a tower of functors to chain complexes, which we denote
$\Pnalg C$ to distinguish from Goodwillie's tower. The functor
$\Pnalg$ approximates
$C$ in the sense that the map $C \rightarrow \Pnalg C$ is the
universal map under $C$ to a functor whose $(n+1)^{\text{st}}$
cross-effect $\Perp_{n+1} C$ is acyclic. In a later work 
\cite{Johnson-McCarthy:deriving-calculus-with-cotriples}, the same
authors show that $\Pnalg$ can be constructed from a cotriple.  One
source of the interest in $\Pnalg$ approximation stems from the fact
that, for many functors, the cross effects are more computationally
accessible than the stabilization; hence, $\Pnalg$ should also be more
accessible than Goodwillie's $P_n$. 

Taking their work as inspiration, we explore to what extent similar
ideas are useful in the study of functors from spaces to spaces. 
In this context, we develop a tower analogous to that of Johnson and
McCarthy. We call this algebraic Goodwillie tower, and denote it
$\Pnalg F(X)$.
The algebraic Goodwillie approximation was created with the intent to
be universal with respect to natural transformations of $F$ to
functors whose $(n+1)^{\text{st}}$ cross effect is contractible; 
however, it turns
out that there is a subtle issue involving $\pi_0$. For example, the
vanishing of the second cross effect forces a monoidal structure on
$\pi_0$, but on $\pi_0$ the approximation process can only produce the
group completion of this monoid (since it involves loops on a space). 
We solve this problem is by requiring our functors
to take values in topological groups or connected spaces.

Our main theorem is that for good functors $F$, 
 there exists a fibration sequence up to homotopy
\begin{equation}
  \label{eq:main-thm-in-intro}
 \sR{ \Perp_{n+1}^{*+1} F } \rightarrow F \rightarrow \Pnalg F ,
\end{equation}
where the fiber is built from a cotriple $\Perp$ formed from the
diagonal of the $(n+1)^{\text{st}}$ cross effect.
In \cite{Johnson-McCarthy:deriving-calculus-with-cotriples},
the authors are able to define the analog of $\Pnalg$ simply by taking
the cofiber of the left-hand map in \eqref{eq:main-thm-in-intro}. 
However, in an unstable topological setting, that approach is not
useful, and the proof is much more difficult. 


It is easy to show that when $F$ commutes with realizations, the
\nth{} algebraic Goodwillie approximation $\Pnalg F(X)$ agrees with the
\nth{} topological Goodwillie approximation $P_n F(X)$ (see
\ref{prop:Pnalg-agrees-Pn}), so in this case
\eqref{eq:main-thm-in-intro} shows that the ``remainder'' of the
$n$-excisive Goodwillie approximation is the cotriple homology of
$F$.


\subsection{Organization of the paper}
We state the main theorem in \S\ref{sec:main-theorem}, 
along with all of the basic definitions necessary to understand it. 
In \S\ref{sec:technical}, we summarize technical
ingredients of the proof, including the limit axiom, $\pi_*$-Kan
functors, results of Goodwillie, and iterated cross effects.
Then \S\ref{sec:Pnalg-connectivity} shows that
the approximation $\Pnalg$ preserves the connectivity of
natural transformations.
To deal with another (important) technical issue,
\S\ref{sec:fiber-contractible-Cartesian} sketches proofs that
when the cross effect vanishes,  a certain cube is actually Cartesian
(rather than just having contractible total fiber).
In \S\ref{sec:perp-is-cotriple}, we prove that $\Perp$ is really a cotriple. 
After these technicalities, we give a detailed outline of the main theorem 
(\S\ref{sec:main-theorem-outline}). The proof of the main
theorem then follows, with \S\ref{sec:perp-F-zero} proving that when
$\Perp_{n+1}F(X)\simeq 0$, then $F(X)\simeq \Pnalg F(X)$, and
\S\ref{sec:perp-F-nonzero} reducing the general case to that case.

\subsection{Acknowledgements}
This paper is based on my Ph. D. thesis at the University of
Illinois. I would like to thank Randy McCarthy 
for his guidance. Thanks also to Jim McClure, Tom Goodwillie, and 
the referee for their helpful suggestions.


\section{Preliminaries}
\label{sec:preliminaries}

The category of ``spaces'' is taken to be the category of compactly
generated spaces with nondegenerate basepoint. 
The space $0$ or $\basept$ is the space
with only one point. When forming the function space $\Map(A,B)$, we
always implicitly replace $A$ by a cofibrant approximation (CW
complex) and $B$ by a fibrant approximation (which is no change with
our definition of spaces). When forming the geometric realization, we
first ``thicken'' each space
(\cite[p.~308]{Segal:categories-and-cohomology-theories}) so the resulting
functor has good homotopy behavior.
When we say two spaces are \emph{equivalent}, we mean they are weakly
homotopy equivalent. 
We use the term \emph{fibration sequence up to homotopy} for a
sequence $F \rightarrow E \rightarrow B$ in which $F
\xrightarrow{\simeq} \hofib(E \rightarrow B)$.

We briefly recall a few standard properties of functors.
A functor is \emph{continuous} if the map $\Map(X,Y)
\rightarrow \Map(F(X),F(Y))$ sending $f$ to $F(f)$ is continuous.
A functor is a \emph{homotopy functor} if it preserves weak homotopy
equivalences. 
A functor is called \emph{reduced} if $F(0)\simeq 0$.
\begin{definition}
\label{def:limit-axiom}
  A homotopy functor $F$ is said to satisfy the \emph{limit axiom} if
  $F$ commutes with filtered homotopy colimits of finite
  complexes. That is, if $\hocolim F(X_\alpha) \xrightarrow{\simeq} F(\hocolim
  X_\alpha)$ for all filtered systems $\Set{X_\alpha}$ of finite complexes.
\end{definition}
In this paper, we assume that all functors are continuous homotopy
functors that satisfy the limit axiom. All functors are defined from
pointed spaces to pointed spaces unless otherwise stated. 

Cubical diagrams and homotopy fibers figure heavily into this work, so
we recall several definitions from \cite{Cal2}.
Let $T$ be a finite set.
$\PP(T)$ is the poset of subsets of $T$
(regarded as a category). $\PP_0(T)$ is the poset of nonempty subsets
of $T$.
A $T$-cube is a functor defined on $\PP(T)$.
If $\cube{X}$ is a $T$-cube of pointed spaces, then its homotopy
fiber, $\hofib \cube{X}$, is the subspace of
the function space
$\prod_{U \subset T} \Map\left([0,1]^U,\cube{X}(U)\right)$
consisting 
of maps that are natural in $U$ and send points with any
coordinate $1$ to the basepoint. 
A formal definition is given in \cite[1.1]{Cal2}.
Alternatively, the homotopy fiber can be constructed by iterating the
process of taking fibers in a single direction.

We write $\mathbf{n}$ for the finite (unpointed) set $\Set{1,\ldots,n}$, 
and $[n]$ for the finite space $\bigvee^n S^0 \cong \Set{0,\ldots,n}$,
with basepoint $0$.


\section{The Main Theorem}
\label{sec:main-theorem}

In this section, we briefly present all of the background necessary to
understand the Main Theorem (\ref{thm:main-theorem}). In brief: cross
effects, left Kan extensions, excisive functors, and $\Pnalg$

The cross effect measures the failure of a functor to take coproducts
to products.
\begin{definition}[$cr_n$] 
Define the $n^{\text{th}}$ cross-effect cube,
$\CRN F(X_1,\ldots,X_n)$, to be the $\mathbf{n}$-cube $\cube{X}$ 
with $\cube{X}(U) = F\left( \bigvee_{u\not\in U} X_u\right)$ and
$\cube{X}(i:U\rightarrow V)$ induced by the identity of $X_u$ if
$u\not\in V$ and the map to the basepoint if $u\in V$. 

Let $cr_n F(X_1,\ldots,X_n)$ denote the homotopy fiber of $\cube{X}$.
\end{definition}
\begin{definition}[$\Perp_n$, $\epsilon$]
Let $\cube{X} = \CRN F(X,\ldots,X)$ and let $\cube{Y}$ be the
$\mathbf{n}$-cube with 
$$\cube{Y}(U) = 
\begin{cases}
F(X) & \text{if $U=\emptyset$} \\
\basept & \text{otherwise}
\end{cases}
$$
Define $\Perp_n F(X) = \hofib \cube{X}$, and note that $F(X) \cong
\hofib \cube{Y}$. Using the fold map $\bigvee X \rightarrow X$ on the
vertex $U=\emptyset$ and the zero
map elsewhere induces a map of cubes $\cube{X} \rightarrow \cube{Y}$. Define 
$\epsilon$ to be the induced map on homotopy fibers, so
$$ \epsilon: \Perp_n F(X) \rightarrow F(X).$$
\end{definition}

As we will show in Section~\ref{sec:perp-is-cotriple}, $\Perp_n$ is
part of a cotriple. Hence there is a simplicial object $Y$ with $Y_k =
\Perp_n^{k+1} F(X)$ and face maps $d_i^{(k)} = \Perp_n^i \epsilon
\Perp_n^{k-i}$. 

\begin{example}
  The second cross effect of the functor $F(X)=Q(X)$ is
  contractible, but the second cross effect of $F(X)=Q(X\wedge X)$ is
  $cr_2 F(X,Y) \simeq Q(X \wedge Y) \times Q(Y \wedge X)$.
\end{example}

The left Kan extension gives a canonical way of extending a functor
from a subcategory to a functor defined on the whole category.
\begin{definition}[left Kan extension]
\label{def:left-Kan}\label{def:left-kan}
Let $LF$ denote the homotopy invariant left Kan extension of a functor
$F$ along the inclusion of finite pointed sets into pointed spaces. Letting
$\cat{S}$ denote a small category of finite pointed sets and $\cat{T}$ denote the
category of pointed spaces, 
the
realization of the following simplicial space can taken to be the
definition of LF(-):
\begin{equation}
  \label{eq:left-kan-spaces}
  [n] \mapsto \bigvee_{(C_0, \ldots, C_n)} F(C_0) \wedge
\left(
    \Map_{\cat{S}}(C_0,C_1) \times \cdots \times \Map_{\cat{T}}(C_n, -
    )
\right)_+ .
\end{equation}
The coproduct is taken over all $(C_0, \ldots, C_n)\in\cat{S}^{\times
  n}$. 
\end{definition}

When $F$ is the restriction of a functor (also called $F$) defined on
$\cat{T}$, then there is a natural map $a: LF(Y) \rightarrow F(Y)$
induced by
$$ 
F(C_0) \wedge \Map(C_0,Y)_+
\rightarrow
F(C_0) \wedge \Map(F(C_0),F(Y))_+
\rightarrow
F(Y)
$$
sending 
$$
y \wedge  f
\mapsto
y \wedge  F(f)
\mapsto
F(f )(y)
.$$
We have required that all functors be continuous to guarantee that the
first map above is continuous. 

\begin{lemma}
When $Y$ is a finite pointed set in $\cat{S}$, the map $a: LF(Y)
\rightarrow F(Y)$ is a simplicial homotopy equivalence.
\end{lemma}
\begin{proof}
When $Y\in\cat{S}$, the category $\cat{S}\downarrow Y$ has a terminal
object. This immediately translates into a homotopy contracting
$LF(Y)$ to $F(Y)\wedge(id_Y)_+ \cong F(Y)$.
\end{proof}

When working with the left Kan extension, we will frequently want to
shift the functor so that we can examine its value on coproducts of
spaces other than $S^0$. To do this, we write $F_X(-)$ for the functor
$F(X\wedge -)$. 

We primarily understand $LF$ as ``$F$ made to commute with
realizations'', in the sense of the following lemma.
\begin{lemma}
\label{lem:LF-commutes-with-realization}
  Let $Y_\cdot$ be a simplicial set. Then 
  $$LF_X(\sR{Y_\cdot}) \simeq \sR{ F_X(Y_\cdot)} = \sR{ F(X \wedge Y_\cdot)}.$$
\end{lemma}
\begin{proof}
  Using \eqref{eq:left-kan-spaces}, this follows from the observation
  that $\Map(S,Y)\cong\prod_{s\in S} Y$ commutes with realizations
  when $S$ is a finite set. 
\end{proof}

Before we define the algebraic Goodwillie tower, we will recall the
classical definition of the topological Goodwillie tower.
\begin{definition}[Cartesian, co-Cartesian]
An $S$-cube $\cube{X}$ is Cartesian if the categorical map
$\cube{X}(\emptyset) \cong \holim_{\PP(S)} \cube{X}
\rightarrow \holim_{\PP_0(S)} \cube{X}$, induced by the inclusion of
the category $\PP_0(S)$ of nonempty subsets of $S$ into the category
$\PP(S)$ of all subset of $S$, is a weak equivalence.
An $S$-cube $\cube{X}$ is co-Cartesian if the categorical map $
\hocolim_{\PP_1(S)} \cube{X} \rightarrow 
\hocolim_{\PP(S)} \cube{X} \cong \cube{X}(S)$, induced by the
inclusion of all proper subsets of $S$ into all subsets of $S$,
is a weak equivalence.
An $S$-cube is strongly co-Cartesian if every sub-$2$-cube is
co-Cartesian. 
\end{definition}
\begin{definition}[$n$-excisive]
(\cite[3.1]{Cal2})
$F$ is $n$-excisive if for every strongly co-Cartesian $(n+1)$-cube
$\cube{X}$, the cube $F\cube{X}$ is Cartesian.
\end{definition}
In \cite{Cal3}, Goodwillie constructs a universal $n$-excisive
approximation $P_n F$ to a functor $F$. The approximations form a tower
of functors equipped with natural transformations of the following form:
$$\xymatrix{ F \ar[r] \ar[rd] & P_n F \ar[d]\\ & P_{n-1} F }$$

We are now in a position to define the algebraic Goodwillie
approximation.
\begin{definition}[$\Pnalg F$]
\label{def:Pnalg-F}
The algebraic Goodwillie approximation $\Pnalg F(X)$ is defined to be
the functor given by applying $P_n$ to the left Kan extension of $F$
shifted over $X$; that is, 
$\Pnalg F(X) = P_n (L F_X) (S^0)$.
\end{definition}
The natural transformation $F(X) \rightarrow \Pnalg F(X)$ arises from
evaluating the map $LF_X \rightarrow P_n (LF_X)$ at $S^0$.

\begin{proposition}
\label{prop:Pnalg-agrees-Pn}
If $F_X$ is a functor that commutes with realizations, that is, the
natural map $LF_X(Y) \rightarrow F_X(Y)$ is an equivalence for 
all spaces $Y$, then 
the natural map $\Pnalg F(X) \rightarrow P_n F(X)$ is an equivalence.
\end{proposition}
\begin{proof}
Given $L F_X \xrightarrow{\simeq} F_X$, applying $P_n$ and evaluating
at $S^0$ gives $P_n(L F_X)(S^0) \xrightarrow{\simeq} P_n F_X(S^0)$. 
The left hand side is 
$\Pnalg F(X)$, and the right hand side is $P_n F(X)$.
\end{proof}

There remains a technical hypothesis on $F$, related to the
``$\pi_*$-Kan condition'', needed in the main theorem. We give
the hypotheses here, but defer discussion of how they are used until
Section~\ref{sec:pi-star-Kan}.  

\begin{hypothesis}[Connected Values]
\label{hyp:connected-values}
\label{hypothesis-1}
$F$ has connected values (on coproducts of $X$) if the functor $L F_X$
is always connected.
\end{hypothesis}

Let $\cat{G}$ denote the category of grouplike
$H$-spaces. Specifically, by $\cat{G}$, we mean the category of
algebras over the associativity operad with inverses and identity up
to homotopy. In this category, all morphisms strictly preserve all
homotopies, so it is rigid enough that the realization of a simplicial
object in $\cat{G}$ is still in $\cat{G}$.

\begin{hypothesis}[Group Values]
\label{hyp:group-values}
\label{hypothesis-2}
In the following definition, let $\cat{T}$ denote the category of
pointed spaces.
Let $U: \cat{G} \rightarrow \cat{T}$ be the forgetful functor.

$F$ ``has group values'' or ``takes values in groups'' if there exists
a functor $F'$ so that the following diagram commutes:
$$\xymatrix{
&
\cat{G}
\ar[d]^U
\\
\cat{T}
\ar[ur]^{F'(-)}
\ar[r]_{F(-)}
&
\cat{T}
}$$
In this case, we will conflate $F( -)$ with its lift to groups. 
We say that $F$ has group values on coproducts of $X$ if $F_X$ has
group values.
\end{hypothesis}

We are now able to state our theorem relating the cross effects of a
functor and the algebraic Goodwillie approximation.
\begin{theorem}[Main Theorem]
\label{thm:main-theorem}
Let $F$ be a reduced homotopy functor from pointed spaces to pointed spaces.
If $F$ has either connected values (\ref{hypothesis-1}) or
group values (\ref{hypothesis-2}) on coproducts of $X$,
then the following is a fibration sequence up to homotopy:
\begin{equation}
\label{eq:perp-fibration}
\realization{\Perp_{n+1}^{*+1} F(X)}
\rightarrow 
F(X)
\rightarrow
\Pnalg F(X)
,
\end{equation}
and furthermore the right hand map is surjective on $\pi_0$.
\end{theorem}

To establish this theorem, we use a ``ladder'' induction on $n$, where
there are two cases for each $n$: the first depends on the second for
smaller $n$, and the second depends on the first for the same $n$.
The first case is to establish that when $\Perp_{n+1} F$
vanishes---that is, when $F$ is degree $n$--- the map $F \rightarrow
\Pnalg F$ is actually an equivalence; this proof proceeds by examining
a fibration sequence obtained from the second case
for smaller $n$. The second case is to consider the general case,
in which $\Perp_{n+1} F$ may be nonvanishing, and proceed by applying
the first case (for the same $n$) to the fiber of the map from $F\rightarrow
\Pnalg F$.
A more extensive outline of the proof is given in
Section~\ref{sec:main-theorem-outline}.


\section{Technical Conditions}
\label{sec:technical}

In this section, we address many technical aspects needed to make our
machinery work.

\subsection{Limit axiom}

It is of primary importance to know that the functors we will be
working with satisfy the limit axiom.
\begin{lemma}
  For any functor $F$,
  the functor $LF$ satisfies the limit axiom (\ref{def:limit-axiom}).
\end{lemma}
\begin{proof}
  The functor $\Map(K,-)$ satisfies the limit axiom for any compact
  $K$. Examining \eqref{eq:left-kan-spaces}, we see that this implies
  that $LF$ satisfies the limit axiom.
\end{proof}

\subsection{Eilenberg-Zilber}

An ``Eilenberg-Zilber''-type theorem for bisimplicial sets provides a
source for connectivity estimates.

\begin{lemma}
\label{lem:basic-Eilenberg-Zilber}
Let $X_\cdot$ and $Y_\cdot$ be simplicial spaces satisfying the $\pi_*$-Kan
condition, and let $f_\cdot: X_\cdot \rightarrow Y_\cdot$ be a map
between them. If $n\ge 0$ and $w \ge -1$ are integers such that for
$i<n$, the map $f_i$ is a weak equivalence, and for $i\ge n$, the map
$f_i$ is $w$-connected, then $\realization{f_\cdot}$ is $(n+w)$-connected.
\end{lemma}
\begin{proof}
Using the homotopy spectral sequence of
\cite[Theorem~B.5]{Bousfield-Friedlander:Gamma-Spaces}, we have a
spectral sequence $E_{*,*} \Rightarrow \pi_*\realization{X_\cdot}$
and a spectral sequence $F_{*,*} \Rightarrow
\pi_*\realization{Y_\cdot}$, with a $f$ inducing a map of spectral
sequences $E \rightarrow F$. The hypotheses imply that $E^1_{i,j}
\cong F^1_{i,j}$ when $i<n$ or $j<w$, and $E^1_{i,w}$ surjects onto
$F^1_{i,w}$. An easy analysis of possible differentials then shows
that the map $E^{\infty}_{i,j} \rightarrow F^{\infty}_{i,j}$ is an
isomorphism when $i+j < n+w$, and a surjection when $i+j = n+w$.
\end{proof}


\begin{corollary}[Eilenberg-Zilber connectivity estimate]
\label{cor:eilenberg-zilber-connectivity-estimate}
  Let $X$ and $Y$ be $\pi_*$-Kan functors 
  from $(\Delta^{\text{op}})^{\times N}$
  to spaces, and let $p=(p_1,\ldots,p_N) \in
  (\Delta^{\text{op}})^{\times N}$ denote an index for these
  multisimplicial spaces.
  Let $f: X\rightarrow Y$.
  Suppose that $f(p)$ is $w$-connected for all indices $p$.
  If in addition, there are integers $\Set{n_i \suchthat i=1,\ldots,N}$ so that
  $f(p)$ is an equivalence if there exists an $i$ with $1\le i \le N$
  such that $p_i < n_i$, then 
  $\sR{f}$ is $(\Sigma n_i + w)$-connected.
\end{corollary}
\begin{proof}
Iterating Lemma~\ref{lem:basic-Eilenberg-Zilber} produces this result.
\end{proof}

\subsection{The $\pi_*$-Kan condition}
\label{sec:pi-star-Kan}

For us to be able to say anything useful about $L F$, we need to
know that the Kan extension of the fiber of a map is the fiber of the
Kan extensions, and that one can continue to say similar things about
the functor defined by taking fibers of certain maps.
\begin{definition}[$\pi_*$-Kan functor]
\label{def:pi-star-Kan-functor}
A functor $F$ is called a $\pi_*$-Kan functor if given a map of
simplicial sets $p: Y_\cdot \rightarrow Z_\cdot$ with a section:
\begin{enumerate}
\item the simplicial spaces $F(Y_\cdot)$ and $F(Z_\cdot)$ satisfy the
$\pi_*$-Kan condition;

\item $\pi_0 F(p_\cdot)$ is a fibration of simplicial sets; and

\item as a functor of $p: Y\rightarrow Z$, the fiber of $F(p)$ is still a
$\pi_*$-Kan functor.
\end{enumerate}
\end{definition}
This is useful in practice due to the following
theorem of Bousfield and Friedlander, restated for simplicial
spaces.

\begin{theorem}
(\cite[Theorem~B.4]{Bousfield-Friedlander:Gamma-Spaces})
\label{thm:bousfield-friedlander}
Let $A$, $B$, $X$, and $Y$ be simplicial spaces, and
suppose that the cube 
$$\xymatrix{
A \ar[r] \ar[d] & X \ar[d] \\
B \ar[r]        & Y }$$
has the property that evaluation at every $[n]\in \Delta^{\text{op}}$
produces a Cartesian cube. If $X$ and $Y$ satisfy the $\pi_*$-Kan condition 
and
if $\pi_0 X \rightarrow \pi_0 Y$ is a fibration of
simplicial sets, then after realization the cube is still Cartesian.
\end{theorem}

In this work we restrict ourselves to 
functors satisfying the hypotheses of connected values
(\ref{hyp:connected-values}) or group values (\ref{hyp:group-values})
so that the simplicial spaces involved always satisfy the 
$\pi_*$-Kan condition.
The $\pi_*$-Kan condition is a technical fibrancy condition introduced
in \cite[\S{}B.3]{Bousfield-Friedlander:Gamma-Spaces} that we do
not repeat here. 

\begin{corollary}
\label{cor:pi-star-fibers-levelwise}
If $F$ is a $\pi_*$-Kan functor, then given a map $p: X\rightarrow Y$
with a section, 
the spaces $L\hofib F(p)$ and $\hofib L F(p)$ are equivalent.
In view of Lemma~\ref{lem:LF-commutes-with-realization}, this is
implies $\sR{ \hofib F(p_\cdot)} \simeq \hofib F (\sR{p_\cdot})$.  
\end{corollary}
\begin{proof}
By definition, a $\pi_*$-Kan functor causes $F(X_\cdot) \rightarrow
F(Y_\cdot)$ to satisfy the hypotheses on the right hand vertical map
in Theorem~\ref{thm:bousfield-friedlander}, so letting $B=\basept$ and
$A = \hofib F(p_\cdot)$ produces the desired result.
\end{proof}

\begin{remark}
  Theorem~\ref{thm:bousfield-friedlander} implies that if each
  space $X_i$ is connected then $\sR{\Omega X_i} \simeq \Omega \sR{X_i}$.
\end{remark}

\begin{lemma}
\label{lem:connected-or-gp-is-pi-star-kan}
  If either $F$ has connected values (\ref{hyp:connected-values})
  or group values (\ref{hyp:group-values}), then
  $F$ is a $\pi_*$-Kan functor (\ref{def:pi-star-Kan-functor}).
\end{lemma}
\begin{proof}
  In these cases, $LF$ always satisfies the $\pi_*$-Kan condition
  (\cite[p.~120]{Bousfield-Friedlander:Gamma-Spaces}). Also,
  for any surjective map $p$, the function $\pi_0 F(p_\cdot)$ is
  a fibration of simplicial sets, since surjective maps of
  simplicial groups are fibrations. If $F$ has connected
  values, then the requirement that $p$ have a section implies, using
  the long exact sequence on homotopy, that
  $\hofib F(p)$ is also a functor to connected spaces. If $F$ has
  group values, then the fact that taking products commutes with
  taking fibers means $\hofib F(p)$ still has group values.
\end{proof}

The cross effect of the Kan extension of a $\pi_*$-Kan functor can be
computed from finite sets, as the following lemma shows.
\begin{lemma}
\label{lem:perp-commutes-with-realization-v2}
Let $Y_1, \ldots, Y_n$ be spaces.
  If $F$ is a $\pi_*$-Kan functor, then 
$$
cr_n (LF)(Y_1,\ldots,Y_n)
\simeq 
L^n (cr_n F) (Y_1,\ldots,Y_n)
,
$$
where $L^n$ indicates the Kan extension is taken in each of the $n$
variables of $cr_n F$.
The statement above can be abbreviated to $\Perp_n(LF)(Y) \simeq L
(\Perp_n F)(Y)$. 
\end{lemma}
\begin{proof}
  All of the maps in $\CRN (LF)$ have sections, so the hypothesis that
  $F$ is a $\pi_*$-Kan functor means that
  Corollary~\ref{cor:pi-star-fibers-levelwise} applies, so taking
  fibers in one direction commutes with realizations. The sections are
  all coherent, so they produce a section on the fibers. 
  The property of being a $\pi_*$-Kan functor passes to the fibers (by
  definition \ref{def:pi-star-Kan-functor}(3)), so the above argument
  applies inductively. 
\end{proof}

\begin{lemma}
\label{lem:perp-commutes-with-realization-for-simplicial-functors}
Let $\Perp = \Perp_n$ for some $n$.
Let $F_\cdot$ be a simplicial object in the category of either
functors to connected spaces or functors to $\cat{G}$ (the category of
``grouplike $H$-spaces''; see definition prior to \ref{hyp:group-values}). 
Then $\realization{\Perp F} \xrightarrow{\simeq} \Perp\realization{F}$.
\end{lemma}
\begin{proof}
As in Lemma~\ref{lem:perp-commutes-with-realization-v2}, the question
is whether the $\Perp$ construction may be taken dimensionwise. By the
same reasoning as Lemma~\ref{lem:connected-or-gp-is-pi-star-kan}, the
spaces in the cube $\CRN F(X, \ldots, X)$ satisfy the $\pi_*$-Kan
condition, and on $\pi_0$ the maps are fibrations of simplicial sets,
so the fibers may be taken dimensionwise. Furthermore, as in that
lemma, the cube of fibers has the same property, so we may proceed
inductively. 
\end{proof}

\subsection{Goodwillie calculus: classification of homogeneous functors}

The functor $P_n$ gives the universal $n$-excisive approximation to a
functor. The functor $D_n$ gives the homogeneous $n$-excisive part of
a functor; it is part of a fibration sequence
$$ D_n \rightarrow P_n \rightarrow P_{n-1} .$$
Goodwillie shows that there is actually a functorial delooping
of the derivative, so this fibration sequence can be delooped once:
\begin{theorem}
\label{thm:delooping-Dn}
(\cite[Lemma~2.2]{Cal3})  
If $F$ is a reduced homotopy functor from spaces to spaces, 
then the map $P_n F \rightarrow P_{n-1} F$ is part of a fibration sequence 
$$ P_n F \rightarrow P_{n-1} F \rightarrow \Omega^{-1} D_n F ,$$
where $\Omega^{-1} D_n F$ is a homogeneous $n$-excisive functor.
\end{theorem} 

\begin{definition}[Derivative of $F$]
\label{def:derivative-of-F}
The $n^{\text{th}}$ derivative of $F$ (at $\basept$), denoted
$\partial^{(n)} F(\basept)$, is the following spectrum with
$\Sigma_n$ action, which we will denote $\mathbf{Y}$.
The space $Y_k$ in the spectrum is
$\Omega^{k (n-1)} cr_n F(S^k,\ldots, S^k)$. The
structure map $Y_k \rightarrow \Omega Y_{k+1}$ arises from suspending
inside and looping outside each variable of $cr_n$. 
\end{definition}

When $F$ satisfies the limit axiom
(\ref{def:limit-axiom}), we can express $D_n F(X)$ using the derivative:

\begin{theorem} (\cite[\S 5, p.~686]{Cal3})
\label{thm:def-of-deriv}
If $F$ is a homotopy functor from spaces to spaces that
satisfies the limit axiom (\ref{def:limit-axiom}),  then
the functor $D_n F$ is given by 
\begin{equation*}
D_n F(X) \simeq \Omega^\infty \left( \partial^{(n)} F(\basept)
\wedge_{h \Sigma_n} X^{\wedge n} \right).
\end{equation*}
where smashing over $h \Sigma_n$ denotes taking homotopy orbits.
\end{theorem}

\subsection{Iterated cross effects}

Classically, cross-effects are defined inductively, by the repeated
application of $cr_2$ in a single variable. With our definition, it is
a theorem that applying $cr_2$ in a single variable of $cr_n$ produces
$cr_{n+1}$.

Recall Goodwillie's notation for sub-cubes: given an $S$-cube
$\cube{W}$ and a subset $A$ of $S$, 
$\partial^A \cube{W}$ denotes the $A$-cube given by restricting
$\cube{W}$ to
subsets of $A$, and 
$\partial_A \cube{W}$ denotes the $(S-A)$-cube with $\partial_A
\cube{W}(B) = \cube{W}(A\cup B)$.

\begin{lemma}
\label{lem:crn-from-cr2}
For $n\ge 2$, there is a natural equivalence
$$cr_2 (cr_n F(X_1, \ldots, X_{n-1}, -))(X_{n},X_{n+1}) 
\xrightarrow{\simeq}
cr_{n+1} F(X_1, \ldots, X_{n+1}).$$
\end{lemma}
\begin{proof}
\begin{raggedright}
Let $S = \mathbf{n}\amalg\mathbf{2}$, and let $T = S -
\Set{1}\amalg\emptyset$. The $S$-cube $\cube{W}$ defining 
$cr_2 (cr_n F(X_1, \ldots, X_{n-1}, -))(X_{n},X_{n+1})$ is
\begin{equation*} 
\cube{W}[X_1,\ldots,X_{n+1}](U\amalg V)
=
    \begin{cases}
    F\left( 
    \bigvee_{v\not\in V} X_v
    \vee
    \bigvee_{u\not\in U\cup\Set{1}} X_{u+1}
    \right)
    &
    \text{$1\not\in U$}
    \\
    F\left( 
    \bigvee_{u\not\in U} X_{u+1}
    \right)
    &
    \text{$1 \in U$}
    \end{cases}
\end{equation*}
\end{raggedright}
Notice that the cube used to compute
$cr_{n+1}F(X_1,\ldots,X_{n+1})$ is exactly the $(n+1)$-cube 
$\partial^T \cube{W}$; that is,
\begin{equation}
\label{eq:a-cube-1}
\hofib \partial^T \cube{W} = cr_{n+1} F(X_1,\ldots,X_{n+1}).
\end{equation}
Also, when $1\in U$, the sub-cube
$\cube{W}(U\amalg -)$ is a constant cube, so 
so
\begin{equation}
\label{eq:a-cube-2}
\hofib \partial_{\Set{1}\amalg \emptyset} \cube{W} \simeq \basept.
\end{equation}

$\cube{W}$ can be written as a
$1$-cube of $(n+1)$-cubes: 
$$\partial^T \cube{W} 
\rightarrow 
\partial_{\Set{1}\amalg \emptyset} \cube{W},$$
so the total homotopy fiber of $\cube{W}$ is the homotopy fiber of the
homotopy fibers of these cubes. Applying \eqref{eq:a-cube-1}
and \eqref{eq:a-cube-2} gives a homotopy fiber sequence
$$cr_2 (cr_n F(-,X_3,\ldots,X_{n+1}))(X_{1},X_{2})
\rightarrow
cr_{n+1} F(X_1,\ldots,X_{n+1})
\rightarrow
\basept,
$$
so the map from the fiber to the total space is a weak equivalence, as
desired.
\end{proof}

\begin{corollary}
\label{cor:cr2-crn-contractible}
Suppose that $\Perp_{n+1} F \simeq \basept$.
Then $cr_2$ applied in any variable of $cr_n F$ results in a
contractible functor.
\end{corollary}
\begin{proof}
First, for any spaces $X_1,\ldots, X_{n+1}$, the space
$cr_{n+1} F(X_1,\ldots, X_{n+1})$ is a retract of $\Perp_{n+1} F(X_1
\vee \cdots \vee X_{n+1}) \simeq \basept$, so $cr_{n+1} F(X_1,
\ldots, X_{n+1})$ is contractible. Lemma~\ref{lem:crn-from-cr2} then
shows that $cr_2$ applied in any variable of $cr_n$ is equivalent to
$cr_{n+1}$, and hence contractible.
\end{proof}


\section{$\Pnalg$ preserves connectivity}
\label{sec:Pnalg-connectivity}

In this section, we establish a property of fundamental importance
when working with $\Pnalg$: the $n$-additive approximation preserves
the connectivity of natural transformations that satisfy some basic
good properties. Actually, we prove the slightly stronger result that
before evaluation at $S^0$, the functor $P_n L(-)_X$ preserves
connectivity. 

\begin{theorem}
\label{thm:Pnd-preserves-connectivity}
Let $F$ and $G$ be reduced $\pi_*$-Kan functors (\ref{def:pi-star-Kan-functor}).
If $\eta: F\rightarrow G$ is a natural transformation that is
$w$-connected, then the natural
transformation $P_{n} ( L \eta )$ is $w$-connected.
\end{theorem}

Once this theorem is established, we have the following immediate corollary.

\begin{corollary}
\label{cor:Pnd-preserves-connectivity}
Let $F$ and $G$ be reduced $\pi_*$-Kan functors.
If $\eta: F\rightarrow G$ is a natural transformation that is
$w$-connected, then the natural
transformation $\Pnalg ( \eta )$ is $w$-connected.
$\qed$
\end{corollary}


\begin{lemma}
\label{lem:deriv-spectra-connectivity}
Let $\eta: F \rightarrow G$ be a natural transformation of $\pi_*$-Kan
functors. If $\eta$ is $w$-connected, then the induced map of
derivative spectra
$\partial^{(n)} LF(\basept) \rightarrow \partial^{(n)} LG(\basept)$ is
$w$-connected. 
\end{lemma}
\begin{proof}
Using the Eilenberg-Zilber connectivity estimate
(\ref{cor:eilenberg-zilber-connectivity-estimate}), 
for all $n\ge 1$, the map 
 $\realization{\Perp_n F(S^k_\cdot)}  \rightarrow 
\realization{\Perp_n G(S^k_\cdot)}$ is $(nk+w)$-connected. The
derivative spectrum $\partial^{(n)}LF(\basept)$ 
then has as its $k^{\text{th}}$ space the space
$\Omega^{k(n-1)} \Perp_n LF(S^k_\cdot)$, which by 
Lemma~\ref{lem:perp-commutes-with-realization-v2}, is equivalent to
$\Omega^{k(n-1)} \realization{\Perp_n F(S^k_\cdot)}$. On these
spaces the map induced by $\eta$ is $(k+w)$-connected, exactly as
required to produce a $w$-connected map 
$\partial^{(n)} LF(\basept) \rightarrow
\partial^{(n)} LG(\basept)$.
\end{proof}

\begin{corollary}
\label{cor:derivatives-connective}
The derivative spectrum $\partial^{(n)} LF(\basept)$ is connective.
\end{corollary}
\begin{proof}
The map $F \rightarrow 0$ is always $0$-connected.
\end{proof}

\begin{corollary}
\label{cor:Dnd-preserves-conn}
If $\eta: F \rightarrow G$ is a $w$-connected map of $\pi_*$-Kan
functors, then $D_n (L \eta)$ is $w$-connected.
\end{corollary}
\begin{proof}
Recall from Theorem~\ref{thm:def-of-deriv}
that for any functor $H$, 
$$
D_n (LH) (Y) 
= \LoopInfty \left(
\partial^{(n)} LH (\basept) \wedge_{h\Sigma_n} (Y)^{\wedge n} \right) . 
$$
Taking homotopy orbits and smashing with a fixed space preserves
connectivity, so this is really a question about the connectivity of
the map $\partial^{(n)} LF (\basept) \rightarrow \partial^{(n)} LG
(\basept)$. The required connectivity was established in
 Lemma~\ref{lem:deriv-spectra-connectivity}.
\end{proof}

\begin{corollary}
\label{cor:Pn+1-Pn-surjective}
If $F$ is a reduced $\pi_*$-Kan functor
then the natural transformation 
$P_{n+1} L F\rightarrow P_n LF$ is surjective in $\pi_0$.
\end{corollary}
\begin{proof}
Theorem~\ref{thm:delooping-Dn} (Goodwillie's delooping of $D_n$) shows
that the fibration 
$$D_{n+1} LF \rightarrow P_{n+1} L F \rightarrow P_n LF$$
deloops to a fibration
\begin{equation}
\label{eq:delooping-of-Dn+1}
P_{n+1} LF \rightarrow P_n LF \rightarrow \Omega^{-1} D_{n+1} LF.
\end{equation}
The delooping of $D_{n+1} LF$ consists of smashing with the
suspension of $\partial^{(n+1)} LF(\basept)$ and taking homotopy orbits. 
By Corollary~\ref{cor:derivatives-connective}, the spectrum
$\partial^{(n+1)} LF(\basept)$ is connective, so its suspension is
$0$-connected.  From the long exact sequence on $\pi_*$, this implies
$P_{n+1} LF \rightarrow P_n LF$ is surjective on $\pi_0$.
\end{proof}

\begin{proof}[Proof of Theorem~\ref{thm:Pnd-preserves-connectivity}]
As above, the spectrum $\partial^{(n+1)} LF(\basept)$ in
Equation~\eqref{eq:delooping-of-Dn+1} is  connective, so the map
on the delooping of $D_{n+1}$ 
induced by $\eta$ is $(w+1)$-connected. The map on the fibers
($P_{n+1} (L \eta)$) is therefore $w$-connected.
\end{proof}

A result that is closely related to
Theorem~\ref{thm:Pnd-preserves-connectivity} is the following, which
says that $\Pnalg$ preserves (good) fibrations that are
surjective on $\pi_0$.
\begin{proposition}
\label{prop:LK-group-connected-base}
Given a space $X$ and functors $A$, $B$, and $C$, suppose 
$A(Y) \rightarrow B(Y) \rightarrow C(Y)$ is a fibration sequence for
all finite coproducts $Y = \bigvee X$ of $X$. 
If, on finite coproducts of $X$, either:
\begin{enumerate}
\item $C$ takes connected values
(Hypothesis~\ref{hypothesis-1});
or
\item $B$ and $C$ take group values (Hypothesis~\ref{hypothesis-2}),
and the map $B\rightarrow C$ is a surjective homomorphism of groups, 
\end{enumerate}
then
\begin{equation*}
LA_X(Z) \rightarrow LB_X(Z) \rightarrow LC_X(Z)
\end{equation*}
is a fibration sequence for all spaces $Z$. Furthermore, the sequence
is surjective on $\pi_0$.
\end{proposition}
\begin{proof}
  This is an easy application of 
  Theorem~\ref{thm:bousfield-friedlander}. 
\end{proof}

\begin{corollary}
  \label{cor:Pnd-group-connected-base}
  Under the conditions of
  Proposition~\ref{prop:LK-group-connected-base}, 
$$\Pnalg A \rightarrow \Pnalg B \rightarrow \Pnalg C$$
is a fibration and the map to the base is surjective on $\pi_0$.
\end{corollary}
\begin{proof}
  By definition, $\Pnalg F(X) = P_n (LF_X)(S^0)$, so 
  combining Proposition~\ref{prop:LK-group-connected-base} with
  Theorem~\ref{thm:Pnd-preserves-connectivity} gives the desired
  result. 
\end{proof}


\section{Fiber contractible implies Cartesian}
\label{sec:fiber-contractible-Cartesian}

In this section, we will sketch proofs of the critical but mainly
technical fact that in
the cases we consider, the cross effect vanishing is equivalent to the
cross effect cubes being Cartesian. We generally want to use the fact
that the cross effect is contractible to conclude that the initial
space in the cross-effect cube is equivalent to the (homotopy) inverse
limit of the rest of the spaces. Unfortunately, this is not always
true; the problem is that the homotopy fiber does not detect failure
to be surjective on $\pi_0$. Here we show that our hypotheses on $F$
are sufficient so that this does not happen.

\begin{lemma}
\label{lem:connected-cartesian}
Let $F$ be a functor satisfying Hypothesis~\ref{hypothesis-1}
(connected values) on coproducts of $X$. 
If $\Perp_n F(X) \simeq 0$, then the cube $\CRN F(X,\ldots,X)$
defining $\Perp_n F(X)$ is Cartesian.
\end{lemma}
\begin{proof}[Proof sketch.]
The first step is to consider the pullback diagram 
$$ X \xrightarrow{p} Z \xleftarrow{q} Y $$
when the spaces $X$, $Y$, and $Z$ are connected, and $p$ and $q$ have
sections. In this case, one can directly show that $\pi_0$ of the
homotopy inverse limit is $0$. 

In general, decompose the desired inverse limit into iterated
pullbacks of the form in the first step. All of the maps have sections
because of the very special form of the cube $\CRN F(X,\ldots,X)$.

This shows that the map from the initial space to the homotopy inverse
limit is (trivially) surjective in $\pi_0$; then Cartesian-ness
follows because the total fiber is contractible.
\end{proof}

\begin{lemma}
\label{lem:groups-cartesian}
Let $F$ be a functor satisfying Hypothesis~\ref{hypothesis-2} (group
values) on coproducts of $X$. 
If $\Perp_n F(X) \simeq 0$, then the cube $\CRN F(X)$ 
defining $\Perp_n F(X)$ is
Cartesian. 
\end{lemma}
\begin{proof}[Proof sketch.]
This strategy here is to decompose $F$ into a fibration involving
$F_0$, the connected component of the identity, and $\pi_0 F$, a
functor to discrete groups.
$$ F_0 \rightarrow F \rightarrow \pi_0 F .$$
The statement for $F_0$ follows from
Lemma~\ref{lem:connected-cartesian}. The statement for $\pi_0 F$ is
Lemma~\ref{lem:groups-surj-X0-lim}, below. Then a short argument shows
no problem arises on $\pi_0$ in $\CRN F(X)$.
\end{proof}

The case of discrete groups is key, so we give a complete proof for
this case. 
\begin{lemma}
\label{lem:groups-surj-X0-lim}
Let $\cube{X}$ be an $n$-cube ($n\ge 1$) 
in the category of discrete groups.
If $\cube{X}$ has compatible 
sections to all structure maps (for example, $\cube{X}=\CRN F(X)$), 
then the map
$$\cube{X}(\emptyset) \rightarrow \lim_{U\in\Power_0(\mathbf{n})} \cube{X}(U)$$
is surjective.
\end{lemma}
\begin{proof}
We need to show that the map above is surjective.
This is equivalent to showing that there exists an
$x_\emptyset \in \cube{X}(\emptyset)$ mapping to each coherent system of
elements $x_U \in \cube{X}(U)$, with $U \not= \emptyset$.

$\cube{X}$ is a cube of groups, and hence all of the structure maps are
group homomorphisms. This allows us to subtract an arbitrary
$w\in\cube{X}(\emptyset)$ from $x_\emptyset$, and subtract the images
$\Image_U(w)$ of $w$ in $\cube{X}(U)$ from each $x_U$, to show the
question is equivalent 
to the existence
of an $x_\emptyset - w \in \cube{X}(\emptyset)$ mapping to
each coherent system of elements $x_U - \Image_U(w)$. 

Given a coherent system of elements $x_U$ in an $n$-cube, let $w$ be
the image of $x_{\Set{n}}$ in $\cube{X}(\emptyset)$ using the section
map $\cube{X}(\Set{n}) \rightarrow \cube{X}(\emptyset)$.
Define $z_U = x_U - \Image_U(w)$, noting that when $\Set{n}\subset V$, we
have $\Image_V(w) = x_W$, so $z_W = 0$.
By the preceding paragraph, the surjectivity that we are trying to establish
is equivalent to the existence of a
$z_\emptyset$ mapping to each coherent collection $z_U$.

If $n=1$, then $\lim_{U\in\Power_0(\mathbf{1})} \cube{X}(U) = X(\Set{1})$, 
so the section map
$\cube{X}(\Set{1}) \rightarrow \cube{X}(\emptyset)$ produces a
$z_\emptyset$ mapping to $z_{\Set{1}}$, as desired.

If $n>1$, then we proceed by induction, assuming the lemma is true for
smaller $n$.
Taking the fiber of $\cube{X}$ in the direction of $\Set{n}$, we have
an $(n-1)$-cube 
$$\cube{Y}(U) := \fib \left(
\cube{X}(U) \rightarrow \cube{X}(U \cup \Set{n}) \right).$$
The cube $\cube{Y}$ satisfies the hypothesis of the lemma because
taking fibers preserves compatible sections.
Notice that for $\Set{n} \not\subset U$, the element $z_U$ passes to
the fiber, since it maps to $z_{U\cup \Set{n}} = 0$. 

Now $\cube{Y}$ is an $(n-1)$-cube, so by induction, the map from
$\cube{Y}(\emptyset)$  to $\lim \cube{Y}(U)$ is surjective. That is,
there exists a $y\in \cube{Y}(\emptyset)$ with $\Image_U(y) = z_U$. 
Mapping $y$ to $z \in \cube{X}(\emptyset)$ gives an element $z$ with
$\Image_U(z) = z_U$ for $U \subset \Set{1,\ldots,n-1}$. As above, if
$\Set{n} \subset U$, then $z_U = 0$, so $\Image_U(z) = z_U$ in this
case as well. Therefore, we have produced an element $z$ mapping to
each coherent collection of elements $z_U$, as desired.
\end{proof}

\section{$\Perp$ is a cotriple}
\label{sec:perp-is-cotriple}

In this section, we verify the claim that $\Perp$ is actually a
cotriple on the category of homotopy functors from pointed spaces to
pointed spaces. 

Let $\cat{S}$ denote the full skeletal category of finite
sets whose objects are $\mathbf{k} = \Set{1,\ldots,k}$, for each
integer $k\ge 0$. 
Given an integer $n$, a homotopy functor
$F$, and a space $X$, we construct a functor
$C$, defined on $\cat{S}^{\text{op}}$ and natural in $F$, $X$, and $n$,
such that $C(\mathbf{k})$ is naturally isomorphic (actually
equivariantly homeomorphic) to  
$\Perp_n^k F(X)$. The functoriality of $C$ shows immediately that
the iterates of $\Perp_n$ assemble to form a simplicial object.

The difference between $\Perp_n^k$ and $C(\mathbf{k})$ is the
difference between computing the homotopy fiber of a cubical diagram
by an iterative process and in a single step. Using a good model
for the total homotopy fiber, these are isomorphic.

Since the strict commutativity of certain diagrams is essential, we
begin by recalling Goodwillie's model for the homotopy fiber of a
cubical diagram.
\begin{definition}[homotopy fiber] (\cite[Definition~1.1]{Cal2})
Let $I=[0,1]$ denote the unit interval. For any set $S$, let $I^S$
denote the product of $\abs{S}$ copies of $I$. For clarity, when
$S=\emptyset$, we take $I^S=\Set{0}$. 

Let $\cube{X}$ be an $S$-cube of spaces. Define the homotopy fiber of
$\cube{X}$ to be the subspace of the space of maps 
$$\hofib \cube{X} \subset \prod_{U\subset S} \Map (I^U,
\cube{X}(U))$$
with a map $\Phi \in \hofib \cube{X}$ satisfying:
\begin{enumerate}
\item $\Phi$ is natural with respect to $S$
\item if $x = (x_1,\ldots,x_u)$ is a point in $I^U$ with some $x_i = 1$,
then  $\Phi_U(x) = \basept$, the basepoint of $\cube{X}(U)$.
\end{enumerate}
\end{definition}

The utility of this definition is that computing the homotopy fiber of
a cube by iterating the process of taking fiber of maps produces a 
result isomorphic (that is, homeomorphic) to the result of taking the
homotopy fiber in one step. 
\begin{lemma}
\label{lem:hofib-homeo-iterated-hofib}
Let $\cube{X}$ and $\cube{Y}$ be cubical diagrams of spaces, and
define the cube of cubes $\cube{Z} = \cube{X}\rightarrow \cube{Y}$.
Then there is a natural homeomorphism 
$$\hofib \cube{Z} 
\xrightarrow{\cong} 
\hofib \left( \hofib \cube{X} \rightarrow \hofib \cube{Y} \right).$$
\end{lemma}
\begin{proof}
Given the above definition of homotopy fiber, this is easy to check;
it follows from the homeomorphism $\Map(I^n,\Map(I^m,X))
\cong \Map(I^{n+m}, X)$.
\end{proof}
This definition of homotopy fiber has one other important property:
\begin{lemma}
Let $\cube{X}$ be an $S$-cube, and let $0$ denote the $S$-cube that is
the one point space for all $U\subset S$. Let $\cube{Y} = \cube{X}\rightarrow
0$ be the $(S\sqcup \Set{\basept})$-cube created using the
zero map to join them. Then there is a homeomorphism
$$\hofib \cube{Y} 
= 
\hofib \left(\cube{X} \rightarrow 0\right) 
\xrightarrow{\cong} \hofib\cube{X}$$
induced by the inclusion of $S$ into $S\sqcup\Set{\basept}$.
\end{lemma}
\begin{proof}
In the definition of homotopy fiber, we see that the inclusion of
$S$ into $S\sqcup\Set{\basept}$ induces a projection from 
$$\prod_{U\subset S\sqcup\Set{\basept}} \Map(I^U,\cube{Y}(U))
\rightarrow
\prod_{U\subset S} \Map(I^U,\cube{X}(U)).
$$
The components of ``$\Phi$'' corresponding to sets containing
$\Set{\basept}$ are all the constant map
to the only available point, $\basept$, so the projection is a
homeomorphism. 
\end{proof}

We now proceed with the proof that $\Perp$ is a cotriple.
We begin by defining a ``diagonal'' to encode the information needed to
construct a cube of coproducts and inclusion and projection maps of
the type used to define the cross effect. 

\begin{definition}[Diagonal]
For any sets $S$ and $U$, and 
given a function $f$ from $S$ to $U$, define the set 
$$B_f = \Set{ (s,u) \suchthat u \not= f(s)} \subset S\times U .$$

Define the ``diagonal'' $\Delta(S,U)$ to be the element of $\Power(S
\times U)$ given by the union of all $B_f$:
$$ \Delta(S,U) = \Set{ B_f \suchthat f \in \Hom(S,U)} .$$
%
\end{definition}
Since $\Delta(S,U)$ is naturally isomorphic to $\Hom(S,U)$ by sending
$B_f$ to $f$, it is obviously a functor in $S$ and $U$.

\begin{definition}[Free cube]
  \label{def:free-cube}
  Given sets $U$ and $S$ and a functor $g$ from the discrete category
  $\Delta(S,U)$ to a 
  pointed category with coproducts (for example pointed spaces or
  cubes of pointed spaces), we define
  $\Free(S,U,g)$ to be
  the $(S \times U)$-cube $\cube{X}$
  with vertices
  $$\cube{X}(A) = \bigvee_{\Set{B \in\Delta(S,U) : A \subset B}} g(B)$$
  Morphisms in $\cube{X}$ are induced by the maps $g(B) \rightarrow
  g(B')$ that are the identity if $B=B'$ and the zero map otherwise.
\end{definition}

We now establish that ``free cubes'' are closed under the pullback
operation.
\begin{lemma}
\label{lem:pullbacks-of-free-cubes}
  Let $m$ and $n$ be sets, let the $(m\times U)$-cube $\cube{X}=
  \Free(m,U,g)$ be a free cube, and 
  let $f: n\rightarrow m$ be a function. The $(n\times U)$-cube
  $\Power(f\times 1)^*
  \cube{X}$ is isomorphic to a free cube $\cube{Y} = \Free(n,U,h)$ with 
  $$h(B) = \bigvee_{B'\in \Delta(f,1)^{-1}(B)} g(B') .$$
\end{lemma}
\begin{proof}
  This is a straightforward argument by expanding the definition of
  $h$ in $\cube{Y}(A)$, combining and interchanging the order of
  quantifiers to turn two coproducts into one, and then verifying that
  the resulting indexing set is the same as the indexing set for 
  $\Power(f\times 1)^*\cube{X}$.
\end{proof}

We are now in a position to identify the relationship $\Perp$ has to
the free cube functor.  For the rest of this section, we fix a functor
$F$ and a space $X$.

Given sets $U$ and $S$, let $c_X$ be the function on $\Delta(S,U)$
that has a constant value $X$. Let $C(S)$ be the contravariant functor
of sets $S$ given by 
$$ C(S) = \hofib F \circ \Free(S,U,c_X) .$$

\begin{lemma}
  $C(S)$ is a contravariant functor of $S$.
\end{lemma}
\begin{proof}
Let $\cube{Y} = C(S)$. Given a function $f: S\rightarrow T$, we can
construct the $(S\times U)$-cube $\Power(f\times 1)^* \cube{Y}$. 
The map in the indexing categories induces 
\cite[\S{}XI.9]{Bousfield-Kan:homotopy-limits-completions-and-localizations}
a map on the homotopy fibers:
$$ \hofib_{\Power(T\times U)} F \cube{Y}  
\rightarrow 
\hofib_{\Power(S\times U)} \Power(f\times 1)^* \cube{Y},$$
so it remains to construct a map 
$$
\hofib_{\Power(S\times U)}
F
\circ
\Power(f \times 1)^* \cube{Y} 
\rightarrow 
\hofib_{\Power(S\times U)}
F \circ
\Free(S,U,c_X) .$$
Recall from Lemma~\ref{lem:pullbacks-of-free-cubes} that $\Power(f
\times 1)^* \cube{Y} $ is a free cube with generating function 
$$h(B) = \bigvee_{\Delta(f,1)^{-1}(B)} X .$$
To specify the desired map between free cubes, it suffices to
specify a natural transformation between their ``generating
functions''; the universal map from $\bigvee X$ to $X$ that is the
identity on each $X$ is a natural choice in this case.
\end{proof}

\begin{lemma}
\label{lem:Ck-Perpk}
  Let $n = \abs{U}$. 
  $C(\mathbf{k})$ is equivariantly homeomorphic to $ \Perp^k_{n} F(X)$.
\end{lemma}
\begin{proof}
  This is a straightforward verification that both are naturally
  homeomorphic (by Lemma~\ref{lem:hofib-homeo-iterated-hofib}) to 
  the homotopy fiber of a $\mathbf{k}\times U$-cube whose vertices are
  $$F \left( \bigvee_{i=1}^{i=k} \bigvee_{v_i \not\in V_i} X \right) .$$
\end{proof}

\begin{lemma}
\label{lem:identify-epsilon}
  Let $n = \abs{U}$. 
  Via the homeomorphism of Lemma~\ref{lem:Ck-Perpk}, 
  the map $C(\mathbf{1}) \rightarrow C(\emptyset)$ induced by $C(i:
  \emptyset \rightarrow \mathbf{1})$ corresponds 
  to the map $\Perp_n F(X)
  \rightarrow F(X)$ induced by $F(\bigvee^n X \rightarrow X)$. 
\end{lemma}
\begin{proof}
  This is also straightforward. Let
  $\cube{X}=\Free(\mathbf{1},U,c_X)$ and observe that 
  the map $\Power(i\times 1)^*$ induces
  the identity on $\cube{X}(\emptyset) = \bigvee^n X$.
\end{proof}

\begin{theorem}
\label{thm:perp-is-cotriple}
  The functor $\Perp$ is a cotriple on the category of homotopy
  functors from pointed spaces to pointed spaces. The map $\epsilon:
  \Perp \rightarrow 1$ is induced by the ``fold map'' $\bigvee X
  \rightarrow X$.
\end{theorem}
\begin{proof}
Lemma~\ref{lem:Ck-Perpk} identifies $\Perp^k$
  and $C(\mathbf{k})$, and Lemma~\ref{lem:identify-epsilon} identifies
  the map $\epsilon$. 
  Then the requisite identities follow from the  applying $C$ to the
  following diagrams:
  $$\xymatrix{
    \Set{1,2} 
    & 
    \Set{1} \ar[l]
    \\
    \Set{2} \ar[u] 
    &
    \emptyset \ar[l] \ar[u]
    }
    \qquad
\xymatrix{
{} 
&
\Set{1}
&
{}
\\
\Set{1}
\ar[ur]^=
\ar[r]^{i_1}
&
\Set{1,2}
\ar[u]
&
\Set{2}
\ar[ul]^=
\ar[l]^{i_2}
}
$$
\end{proof}


\section{Main Theorem: Outline}
\label{sec:main-theorem-outline}

To establish the main theorem, we use induction on $n$, beginning with
the case $n=1$.  We further break down the induction into ``Case I'',
where $F$ is degree $n$ --- that is, $\Perp_{n+1} F(X)
\simeq 0$, and ``Case II'', which shows that Case~I implies the result
for arbitrary $F$. 
Our proof will be a ladder induction, with Case~I
depending on Case~II for smaller values of $n$, and Case~II depending
on Case~I for the same value of $n$.

In Section~\ref{sec:perp-F-zero}, we treat the case when
$\Perp_{n+1} F(X) \simeq 0$. In this case, we show directly that the
fiber of the fibration sequence we obtain from induction,
$$
\realization{\psp[n]{F}(X)}
\rightarrow 
F(X)
\rightarrow
P_{n-1}^{\alg} F(X)
,
$$
is a homogeneous degree $n$ functor. This implies that $F(X) \simeq
\Pnalg F(X)$ in this case.

\begin{definition}[$A_F$]
\label{def:AF}
Define the functor $A_F(X)$ to be the homotopy fiber in the
fibration:
\begin{equation}
  \label{eq:AF-def}
A_F(X)
\rightarrow
\rpspfx{F}{X}
\rightarrow 
F(X).
\end{equation}
\end{definition}

For the general case, 
in Section~\ref{sec:perp-F-nonzero}, we consider the
auxiliary diagram: 
$$\xymatrix{
A_F(X)
\ar[r]
\ar[d]
&
\realization{\psp{F}(X)}
\ar[r]^{\epsilon}
\ar[d]
&
F(X)
\ar[d]
\\
\Pnalg A_F(X)
\ar[r]
&
\Pnalg \left( \sR{\psp{F}(X)} \right)
\ar[r]
&
\Pnalg F(X)
}$$
where the bottom row is shown to be a fibration sequence up to homotopy
(Proposition~\ref{prop:perp-f-nonzero-fib-conn}
and~\ref{prop:perp-f-nonzero-fib-discrete}).  
We show that $\Perp_{n+1} A_F(X)
\simeq 0$ (Lemma~\ref{lem:perp-AF-zero}), 
and hence Case~I shows that
there is an equivalence of the fibers, so the square on the right is
Cartesian. Then it is not hard to establish that 
$\Pnalg \left( \realization{\psp{F}(X)}\right) \simeq 0$
(Lemma~\ref{lem:Pnd-perp-contractible}),
so that the sequence in the main theorem
is actually a fibration sequence up
to homotopy.


\section{Main Theorem, Case I: $F$ degree $n$}
\label{sec:perp-F-zero}

In this section, the goal is to establish that when the $(n+1)$-st
cross effect of $F$ vanishes, $F$ is equivalent to its $n$-additive
approximation, $\Pnalg F$. 
The proof of this result
(Proposition~\ref{prop:perp-F-zero-realizations}) will  
be given in Section~\ref{sec:proof-of-perp-F-zero}.

\subsection{Additivization and the bar construction}

The case $n=1$ is the work of Segal.

\begin{lemma}
\label{lem:perp-2-zero-Segal}
Suppose $F$ is a reduced functor that has either connected values
(\ref{hypothesis-1}) or 
group values (\ref{hypothesis-2}) and commutes with realizations. 
If $\Perp_{2} F \simeq 0$, then $F \simeq \Omega \circ F \circ \Sigma$.
\end{lemma}
\begin{proof}
Let $X$ be a space.
Under these hypotheses, the map 
$F(X\vee X) \xrightarrow{\simeq} F(X) \times F(X)$ is an equivalence,
so we can regard $[n] \mapsto F(\bigvee^n X)$ as a $\Gamma$-space.
Segal's work \cite[Proposition~1.4]{Segal:categories-and-cohomology-theories}
then 
shows that $F(X) \simeq \Omega B F(X)$, where $B$ denotes the bar
construction, and note that $BF(X) \simeq \realization{F(X \wedge
  S^1_\cdot)} \simeq F(S^1 \wedge X)$, because $F$ commutes with
realizations. 
\end{proof}

\begin{corollary}
\label{cor:perp-f-spectrum}
If $F$ has either connected values
(\ref{hypothesis-1}) or 
group values (\ref{hypothesis-2}), and
$\Perp_{n+1} F \simeq 0$, then, as a symmetric functor of $n$
variables, $\Perp_n F$ is the infinite loop space of a
symmetric functor to connective spectra $\mathbf{\BoldPerp_n F}$: 
$$\Perp_n F \simeq \LoopInfty \left( \mathbf{\BoldPerp_n
    F}  \right).$$
\end{corollary}
\begin{proof}
When $F$ satisfies the hypotheses above, then $\Perp_{n} F$ also
satisfies these hypotheses. Since $\Perp_{n+1} F\simeq 0$, we know that
applying $\Perp_2$ in any input to $cr_n F$ is
contractible (by Corollary~\ref{cor:cr2-crn-contractible}), 
so Lemma~\ref{lem:perp-2-zero-Segal} gives $\Perp_n
F$ as the first space of a connective $\Omega$-spectrum. 
\end{proof}

\begin{corollary}
\label{cor:perp-2-zero-F-P1dF}
Suppose $F$ is a reduced functor that commutes with realizations and
has either connected values
(\ref{hypothesis-1}) or 
group values (\ref{hypothesis-2}), and
$\Perp_{2} F \simeq 0$, then $F \simeq P_1 F$.
\end{corollary}
\begin{proof}
We need to show that $F \rightarrow \Omega^n \circ F  \circ \Sigma^n$ 
is an equivalence for all $n$. This almost follows from
Lemma~\ref{lem:perp-2-zero-Segal}, but one needs to check that $\Omega
\circ F \circ \Sigma$ commutes with realizations. 
This follows from Theorem~\ref{thm:bousfield-friedlander} 
because $F(S^1\wedge -)$ is always connected.
\end{proof}

\subsection{The equivariant structure of the cross-effects of a
  symmetric functor}
We begin with a definition of the ``symmetric group'' and its action
on the set of integers $\mathbf{n}= \Set{1,\ldots,n}$.
 By $\Sigma_n$, we mean the set of bijective
set maps from $\mathbf{n}$ to $\mathbf{n}$ with the
group operation corresponding to composition of functions.
So no confusion arises, let us fix a
representation of the symmetric group on the set $\mathbf{n}$: let
$\sigma\in\Sigma_n$ act on $\mathbf{n}$ by sending $j$ to
$\sigma(j)$. To verify that this is an ``action'', note that
$\tau_* \sigma_*$ sends $j$ to $\tau(\sigma(j)) = (\tau\circ \sigma)(j)$,
which is exactly $(\tau\sigma)_*$.
Write $\Sigma_n^+$ for the space $\Sigma_n$ with a disjoint basepoint
added. 

\begin{definition}{\cite[p.~675]{Cal3}}
A functor of $n$ variables $H(X_1,\ldots,X_n)$ is called
\emph{symmetric} if for each permutation $\sigma\in\Sigma_n$ there is
a natural transformation 
$$\sigma_*: H(X_1, \ldots, X_n) \rightarrow H(X_{\sigma^{-1}(1)}, \ldots,
X_{\sigma^{-1}(n)})$$
satisfying $(\sigma \tau)_*=\sigma_* \tau_*$. This definition differs
from Goodwillie's only in that we use $\sigma^{-1}$ where he uses
$\sigma$, and consequently the action of a composition appears in a
different order than his.
\end{definition}

We immediately begin using some abbreviated notation.
\begin{notation}[{$H[\sigma]$}]
Given a symmetric functor of $n$ variables $H$, a map of sets
$\sigma: \mathbf{n}\rightarrow\mathbf{n}$, 
and $n$ ambient variables $X_1, \ldots, X_n$, define
$$ H[\sigma] = H(X_{\sigma(1)}, \ldots, X_{\sigma(n)}).$$
\end{notation}

With this notation, the symmetric structure is $\sigma_*: H[1]
\rightarrow H[\sigma^{-1}]$. Notice that $\tau_*: H[\sigma^{-1}]
\rightarrow H[\sigma^{-1}\tau^{-1}] = H[(\tau\sigma)^{-1}]$.

\begin{definition}[$r_\sigma$]
When $X_1 = \cdots = X_n$, let $r_\sigma$ be the ``relabeling'' map
induced by the isomorphisms $X_{i} \mapsto X_{\sigma(i)}$.
\end{definition}

\begin{definition}[Action of $\Sigma_n$ on the diagonal of a symmetric
functor]
The diagonal of a symmetric functor of $n$ variables has a natural
action of the symmetric group on $n$ letters. Let $X = X_1 =\cdots = X_n$,
and let $\sigma\in\Sigma_n$. The action of $\sigma$ on $\diag H(X)$ is
the composite:
$$ 
H[1] 
\xrightarrow{\sigma_*} 
H[\sigma^{-1}] 
\xrightarrow{r_\sigma}
H[1] 
$$
where $r_\sigma$ is the isomorphism $X_i \cong X_{\sigma(i)}$ that is
available because all of the $X_i$ are equal.
\end{definition}

\begin{remark}
\label{remark:action-on-Hsigma}
By relabeling variables, one sees immediately that the action of
$\tau$ on $H[\sigma]$ is given by the composite:
$$H[\sigma] \xrightarrow{\tau_*} H[\sigma\tau^{-1}] 
\xrightarrow{r_{\sigma\tau\sigma^{-1}}} H[\sigma].$$
\end{remark}

\begin{example}
\label{example:coproduct}
The fundamental example is the wedge. When $n=3$, let $H(X_1,X_2,X_3) =
X_1 \vee X_2 \vee X_3$. 
The symmetric structure map for the cycle $\sigma=(1 2 3)$ in
$\Sigma_3$ is
$$\sigma_*: x_1\vee x_2 \vee x_3 \in H(X_1,X_2,X_3) \mapsto
x_3 \vee x_1 \vee x_2 \in H(X_3, X_1, X_2).$$ 
Notice that
$H(X_3,X_1,X_2) = H[\sigma^{-1}]$, and not $H[\sigma]$.
The map $r_\sigma$ then has us regard
$x_3$ as an element of $X_1$, \emph{etc}.
\end{example}

For the rest of this section, we have will work with the diagonals of
symmetric functors, so set $X_1=\cdots=X_n$, but label them
differently to be able to see the action of the symmetric group more
clearly. To keep notation under control, we write $H$ instead of
$\diag H$.

\begin{lemma}
\label{lem:cross-effect-1}
Let $H$ be a symmetric functor of $n$ variables. If in each variable
$i$, the categorical map
\begin{equation}
\label{eq:H-coprod-prod}
H(\ldots, X_i \vee Y_i, \ldots) 
\xrightarrow{\simeq}
H(\ldots, X_i, \ldots) \times H(\ldots, Y_i, \ldots)
\end{equation}
is an equivalence, then
$$\Perp_n H[1] \xrightarrow{\simeq} \prod_{\sigma\in\Sigma_n} H[\sigma].$$
Furthermore, with the following action, the map is equivariant.

Let $\beta\times\alpha\in\Sigma_n\times\Sigma_n$ act on the left as
follows: $\alpha$ acts on $H$, and hence $\alpha$ acts on $\Perp_n H$
(via $\Perp_n \alpha$) because $\Perp_n$
is a functor; $\beta$ acts on $\Perp_n$ because $\Perp_n$ is the
diagonal of the symmetric functor $cr_n$.

Let $\beta\times\alpha$ act on the right as follows: $\alpha$ acts diagonally
on all copies of $H$; $\beta$ acts by permuting coordinates so that
the $\beta\sigma$ component of $\beta_*(h)$ equals the $\sigma$
component of $h$.
\end{lemma}
\begin{proof}
First note that the actions on the left commute because the action on
$\Perp_n$ is a natural transformation. The actions on the right
obviously commute. 

Using the hypothesis in each variable of $H$, we see that
the categorical map to product (by iterating
\eqref{eq:H-coprod-prod}) is an equivalence. 
$$
\myH\left(\bigvee_{i=1}^{n} X_i, \ldots, \bigvee_{i=1}^{n} X_i \right)
\xrightarrow{\simeq}
\prod_{\sigma: \mathbf{n}\rightarrow \mathbf{n}} 
        \myH[\sigma] 
$$
Using this to compute the cross effect, we see there is an equivalence
\begin{equation}
\label{eq:perp-n-H-decompose-1}
\Perp_n \mathbf{H}[1] = \Perp_n \mathbf{H} (X_1, \ldots, X_n)
\xrightarrow{\simeq}
\prod_{\sigma\in\Sigma_n} \myH[\sigma] .
\end{equation}

Also, notice that substituting $X_{\beta^{-1}(i)}$ for $X_i$ gives:
$$
\Perp_n \mathbf{H}[\beta^{-1}] 
\xrightarrow{\simeq}
\prod_{\sigma\in\Sigma_n} \myH[\beta^{-1}\sigma] .
$$

A short reflection on the origins of the map in 
\eqref{eq:perp-n-H-decompose-1} should convince the
reader that it is equivariant with respect to the action of $\Sigma_n$
on $H$.

The equivariance of the action of $\Sigma_n$ on $\Perp_n$ is
a bit more complicated, so we spell it out in detail.
To orient the reader, consider the commutative diagram whose vertical
composite is action by $\beta$ on $\Perp_n H(1)$.
$$\xymatrix{
\Perp_n  H[1]
\ar[r]^{\simeq}
\ar[d]^{\beta_*}
& 
\prod_{\sigma\in\Sigma_n}  H[\sigma] 
\ar@{-->}[d]
\\
\Perp_n  H[\beta^{-1}]
\ar[r]^{\simeq}
\ar[d]^{r_\beta}
&
\prod_{\sigma\in\Sigma_n} H[\beta^{-1}\sigma]
\ar@{-->}[d]
\\
\Perp_n  H[1] 
\ar[r]^{\simeq}
&
\prod_{\sigma\in\Sigma_n}  H[\sigma]
}$$
Under $\beta_*$, the copy of $H[\sigma]$ in the $\sigma$ component of
the top row is
sent to $H[\beta^{-1}\beta\sigma]$ in the $\beta\sigma$ component of
the middle row. This is often a confusing point, but it is made more
evident by the fact
that the map $\beta_*$ exists even if the inputs $X_i$ are all
different, and $\beta_*$ does not even require that $H$ be a
symmetric functor (only that it be additive in each variable).
Under the ``relabeling'' isomorphism $r_\beta$, this becomes
$H[\beta\sigma]$ in the $\beta\sigma$ component of the bottom row.
This shows that the action of $\beta\in\Sigma_n$ on $\Perp_n  H[1]$
corresponds to the action sending the $\sigma$ component of the
product isomorphically to the $\beta\sigma$ component.
\end{proof}

\begin{lemma}
\label{lem:cross-effect-2}
Let $\myH$ be a symmetric functor of $n$ variables taking values in
spectra. There is a stable equivalence
$$ \Sigma_n^+ \wedge \myH[1] \xrightarrow{\simeq}
\prod_{\sigma\in\Sigma_n} \myH[\sigma]$$
given by
$$\Sigma_n^+ \wedge \myH[1]
\cong
\bigvee_{\sigma\in\Sigma_n} \myH[1]
\xrightarrow{\vee r_\sigma}
\bigvee_{\sigma\in\Sigma_n} \myH[\sigma]
\hookrightarrow^{\simeq}
\prod_{\sigma\in\Sigma_n} \myH[\sigma] .
$$
Let $\beta\times\alpha\in\Sigma_n\times\Sigma_n$ act on $\Sigma_n^+
\wedge \myH[1]$ as follows: $\alpha$ acts on $\myH[1]$, and $\beta$
sends $\sigma\wedge h$ to $(\beta\sigma) \wedge h$.

With the action $\Sigma_n\times\Sigma_n$ on the right as in
Lemma~\ref{lem:cross-effect-1}, this map is equivariant.
\end{lemma}
\begin{proof}
The equivariance with respect to the action of $\beta$ is
immediate. The equivariance with respect to the action of $\alpha$
follows from the commutativity of both of the following squares, where
the vertical composites are the action of $\alpha$ on the left and one
component of the right.

$$\xymatrix{
\myH[1]
\ar[r]^{r_\sigma}
\ar[d]^{\alpha_*}
&
\myH[\sigma] 
\ar[d]^{\alpha_*}
\\
\myH[\alpha^{-1}]
\ar[r]^{r_\sigma}
\ar[d]^{r_\alpha}
&
\myH[\sigma\alpha^{-1}] 
\ar[d]^{r_{\sigma\alpha\sigma^{-1}}}
\\
\myH[1]
\ar[r]^{r_\sigma}
&
\myH[\sigma]
}
$$
The map $r_{\sigma\alpha\sigma^{-1}}$ is used for the action of
$\alpha$ on $H[\sigma]$, as noted in
Remark~\ref{remark:action-on-Hsigma}.
The stable equivalence of
coproducts and products is the last map used.
\end{proof}

\begin{corollary}
\label{cor:perp-n-decomp-to-sigma-n}
Let $\myH$ be a symmetric functor of $n$ variables taking values in
spectra. If the second cross effect in each variable of $\myH$ is
contractible, then 
there is an weak equivariant map
\begin{equation*}
\Perp_n \mathbf{H}[1] 
\xrightarrow{\simeq}
\prod_{\sigma\in\Sigma_n} 
\mathbf{H}[\sigma]
\xleftarrow{\simeq}
\Sigma_n^+ \wedge \mathbf{H}[1] 
.
\end{equation*}
The action of $\beta\in\Sigma_n$ on $\Perp_n$ corresponds to
multiplication by $\beta$ on $\Sigma_n^+$ on the right side of the
equation, and the 
action on $\myH$ is the same on both sides.
\end{corollary}
\begin{proof}
Since $\myH$ is a functor to spectra, the vanishing of the second
cross effect implies that the hypotheses of
Lemma~\ref{lem:cross-effect-1} are satisfied.
Combining Lemmas~\ref{lem:cross-effect-1} and~\ref{lem:cross-effect-2}
gives the desired result.
\end{proof}
This type of identification goes
back at least to Eilenberg and MacLane in the 1950s 
\cite[Theorem~9.1]{Eilenberg-MacLane:H-Pi-n-II}.
With this model, the map $\epsilon: \Perp_n \mathbf{H} \rightarrow
\mathbf{H}$ is given by $\sigma \wedge x \mapsto x$.

\begin{corollary}
\label{cor:decomp-identify-epsilon}
Let
 $\overbar{\epsilon}$ be the map $\Sigma_n^+ \wedge \myH[1]
\rightarrow \myH[1]$ sending $\sigma \wedge h$ to $h$.
Under the conditions of Corollary~\ref{cor:perp-n-decomp-to-sigma-n},
the following diagram commutes:
$$\xymatrix{
\Perp_n \myH[1] 
\ar[r]^{\simeq}
\ar[d]^{\epsilon}
&
\prod_{\sigma\in\Sigma_n} \myH[\sigma]
&
\Sigma_n^+\wedge \myH[1]
\ar[l]^{\simeq}
\ar[d]^{\overbar{\epsilon}}
\\
\myH[1]
\ar[rr]^{=}
&
&
\myH[1]
} .$$
\end{corollary}
\begin{proof}
Let $X=X_1=\cdots=X_n$ be spaces.
The map $\epsilon$ is induced by the fold map $\bigvee X_i \rightarrow
X$, so it sends 
$$\myH[\sigma] = \myH(X_{\sigma(1)}, \ldots, X_{\sigma(n)})
\rightarrow
\myH(X,\ldots,X) = \myH(X_1,\ldots,X_n) = \myH[1]$$
by the isomorphism $X_{\sigma(i)} = X = X_i$; this is the map we have
denoted $r_{\sigma^{-1}}$. On $\sigma\wedge h$, the composite of the
upper and left maps is then $r_{\sigma^{-1}} r_\sigma(h) = h$, and
$\overbar{\epsilon}(\sigma \wedge h) = h$, so the diagram commutes.
\end{proof}


\subsection{Iterated cross effects produce homogeneous functors}

\begin{lemma}
\label{lem:perp-homotopy-orbits-spectrum}
If $F$ has either connected values
(\ref{hypothesis-1}) or 
group values (\ref{hypothesis-2}) on coproducts of $X$, and
$\Perp_{n+1} F(X) \simeq 0$, then
$$\realization{\psP[n]{F}(X)}
\simeq 
 \mathbf{\BoldPerp_n F}(X) \wedge_{\Sigma_n} E \Sigma_n^+
,
$$
where $\mathbf{\BoldPerp_n F}$ denotes the lift to spectra of $\Perp_n F$,
as in Corollary~\ref{cor:perp-f-spectrum}.
\end{lemma}

With these hypotheses,
Corollary~\ref{cor:perp-f-spectrum} shows that $\Perp_n F(X) \simeq
\LoopInfty \mathbf{\BoldPerp_n F}(X)$, so we are entitled to consider the
functor to spectra ``$\mathbf{\BoldPerp_n F}$''. 
Corollary~\ref{cor:cr2-crn-contractible} shows that the second cross
effect in each variable of $\mathbf{\BoldPerp_n F}$ is contractible, so we
may apply Corollary~\ref{cor:perp-n-decomp-to-sigma-n} with
$\myH = \mathbf{\BoldPerp_n F}$ to conclude that
\begin{equation}
\label{eq:perp-is-sigma-n}
 \Perp_n \mathbf{\BoldPerp_n F}(X) 
\simeq 
\Sigma_n^+ \wedge \mathbf{\BoldPerp_n F}(X) .
\end{equation}

We are now in a position to understand $\Perp_n^* \mathbf{\BoldPerp_n F}$.
The issue of how multiplication by $\sigma$ arises from the $\epsilon$
above is somewhat subtle, so we spell it out in detail.
Applying~\eqref{eq:perp-is-sigma-n} repeatedly at each level, we have 
$$ \Perp_n^k \mathbf{\BoldPerp_n F}(X)
\simeq
\underbrace{
\Sigma_n^+ \wedge \cdots \wedge \Sigma_n^+
}_{\text{$k$ factors}} \wedge 
\mathbf{\BoldPerp_n F}(X)
.$$
Recall that the face maps from dimension $n$ to $n-1$ are given by
$d_i = \Perp_n^i \epsilon \Perp_n^{n-i}$. In dimension $k$, the face map
$d_k = \epsilon \Perp_n^k$ just drops the first element 
(by Corollary~\ref{cor:decomp-identify-epsilon}):
$$ d_k (g_k \wedge \cdots \wedge g_1 \wedge y) = g_{k-1} \wedge \cdots
\wedge g_1 \wedge y .$$
To compute the others, note that for any $f$, the map $\Perp_n(f)$ is
equivariant with respect to the action of $\Sigma_n$ on $\Perp_n$ (by
permuting inputs), so in particular $\Perp_n(\epsilon): \Perp_n(\Perp_n F)
\rightarrow \Perp_n F$ is equivariant with respect to the action on of
$\Sigma_n$ on the leftmost $\Perp_n$, so 
\begin{align*}
  \Perp_n \epsilon (g \wedge y ) 
&=
\Perp_n\epsilon(g * (1 \wedge y) )
\\
&= g * \Perp_n\epsilon(1 \wedge y) 
\\
&= g * y,
\end{align*}
where the last follows since the degeneracy $\delta: \Perp_n F
\rightarrow \Perp_n^2 F$ which has $\delta(y) = 1\wedge y$ is a section
to the face map $\Perp_n\epsilon$. This argument shows that all of the
face maps $d_j$ with $0\le j<k$ are given by multiplying $g_{j+1}$ by
the next coordinate to the right (either $g_j$ if $j>0$ or $y$ if
$j=0$).

This is a standard model for 
$E \Sigma_n^+ \wedge_{\Sigma_n} \mathbf{\BoldPerp_n F}(X)$, 
so we have proven Lemma~\ref{lem:perp-homotopy-orbits-spectrum}. $\qed$

\subsection{Proof of Main Theorem, Case~I}
\label{sec:proof-of-perp-F-zero}

\begin{proposition}
\label{prop:perp-F-zero-realizations}
Suppose that $F$ commutes with realizations and has connected values
(\ref{hyp:connected-values}) or group values (\ref{hyp:group-values}). 
If $\Perp_{n+1} F \simeq 0$, then $F \simeq P_n F$.
\end{proposition}

\begin{proof}
When $n=0$, the hypothesis that $F$ is reduced makes the result
trivial. 
When $n=1$, Corollary~\ref{cor:perp-2-zero-F-P1dF} shows
that $F(X) \simeq P_1 F(X)$, so that establishes the truth of the base case
in our induction.

Finally, when $n>1$ we apply Proposition~\ref{prop:perp-F-nonzero}
with one smaller $n$ to produce a fibration sequence:
\begin{equation}
\label{eq:fib-seq-with-n-1}
\realization{\Perp_n^* (\Perp_n F) }
\rightarrow
F
\rightarrow
P_{n-1} F
,
\end{equation}
where the map $F\rightarrow P_{n-1} F$ is surjective on
$\pi_0$.  (Recall that since $F$ commutes with realizations, 
$P_{n-1}^{\alg} F \simeq P_{n-1} F$.)

We now show that the fiber here is an $n$-excisive functor. 
From Corollary~\ref{cor:perp-f-spectrum}, we know that $\Perp_n F$
lifts to a functor to connective spectra $\mathbf{\BoldPerp_n F}$, so that 
$\Perp_n F \simeq \LoopInfty \mathbf{\BoldPerp_n F}$. Using this, we have:
\begin{gather*}
\sR{\Perp_n^* (\Perp_n F) } 
\simeq
\sR{ \Perp_n^* \left( \LoopInfty \mathbf{\BoldPerp_n F} \right) }
.
\\
\intertext{The functor $\LoopInfty$ is a right adjoint, so it
  preserves homotopy fibers, and hence commutes with $\Perp_n$:}
\sR{ \Perp_n^* \left( \LoopInfty \mathbf{\BoldPerp_n F} \right) }
\simeq
\sR{  \LoopInfty  \Perp_n^* \left(\mathbf{\BoldPerp_n F} \right) }
.
\\
\intertext{Now $\mathbf{\BoldPerp_n F}$ is a functor to connective
  spectra, and hence all applications of $\Perp_n$ to it result in
  functors to
  connective spectra, so we can use
  \cite{Beck:classifying-spaces-for-homotopy-everything-H-spaces}
  to move $\LoopInfty$ outside of the realization:}
\sR{  \LoopInfty  \Perp_n^* \left(\mathbf{\BoldPerp_n F} \right) }
\simeq
\LoopInfty  \sR{  \Perp_n^* \left(\mathbf{\BoldPerp_n F} \right) }
.
\\
\intertext{Then Lemma~\ref{lem:perp-homotopy-orbits-spectrum} shows
  that the term inside the realization computes the homotopy orbits of
the $\Sigma_n$ action on $\mathbf{\BoldPerp_n F}$:}
\sR{ \Perp_n^* (\mathbf{\BoldPerp_n F}) } 
\simeq
\mathbf{\BoldPerp_n F } \wedge_{\Sigma_n} E \Sigma_n^+
.
\\
\intertext{Combining all of these, we identify the fiber in
  \eqref{eq:fib-seq-with-n-1} as the infinite loop space of the
  preceding homotopy orbit spectrum:} 
\sR{\Perp_n^* (\Perp_n F) } 
\simeq
\LoopInfty \left(
\mathbf{\BoldPerp_n F } \wedge_{\Sigma_n} E \Sigma_n^+
\right)
.
\end{gather*}

We now establish that this functor is $n$-excisive. Since it is known
to be the fiber of $F \rightarrow P_{n-1} F$, this will imply that it
is actually homogeneous $n$-excisive.

The functor $\LoopInfty$ preserves Cartesian squares, so we need only
establish that $\mathbf{\BoldPerp_n F } \wedge_{\Sigma_n} E
\Sigma_n^+$ is $n$-excisive. 
By hypothesis, $\Perp_{n+1} F \simeq \basept$, so
Corollary~\ref{cor:cr2-crn-contractible} shows that 
the second cross effect in any variable of $\mathbf{cr_n F}$
is contractible. Hence, by Corollary~\ref{cor:perp-2-zero-F-P1dF}, in
each variable  $\mathbf{cr_n F}$ is $1$-excisive. Then
\cite[Proposition~3.4]{Cal2} shows that its diagonal,
$\mathbf{\BoldPerp_n F}$ is $n$-excisive. 
That is, given any strongly
co-Cartesian $(n+1)$-cube $\cube{X}$, the map
$$\mathbf{\BoldPerp_n F} \cube{X}(\emptyset) 
\rightarrow
\holim_{\Power_0(\mathbf{n+1})} \mathbf{\BoldPerp_n F} \cube{X}$$
is an equivalence. Since $\Sigma_n$ acts naturally on
$\mathbf{\BoldPerp_n F}$, this is a $\Sigma_n$-equivariant map, so it is
still an equivalence after taking homotopy orbits. This establishes
that the functor
$\mathbf{\BoldPerp_n F } \wedge_{\Sigma_n} E \Sigma_n^+$ is
$n$-excisive, as desired.

Equation~\eqref{eq:fib-seq-with-n-1} is a fibration sequence up to homotopy, 
so we know the natural map 
$$\realization{\Perp_{n}^{*} (\Perp_n F) } 
\rightarrow
\fib ( F \rightarrow P_{n-1} F )
$$
is an equivalence.
Applying $P_n$ and using the fact that we have just shown that
$\realization{\Perp_{n}^{*}
  (\Perp_n F) } $ is $n$-excisive, we have
\begin{align*}
 \realization{\Perp_{n}^{*} (\Perp_n F) } 
&\simeq
 P_n \realization{\Perp_{n}^{*} (\Perp_n F) } 
\\
&\simeq
P_n \fib( F \rightarrow P_{n-1} F ) 
\\
&\simeq 
\fib ( P_n F \rightarrow P_{n-1} F )
\\
&\simeq 
 D_n F .
\end{align*}
When $F$ commutes with realizations,
the map $F \rightarrow P_{n-1} F$ is 
surjective on $\pi_0$ (Theorem~\ref{thm:Pnd-preserves-connectivity}).
This lets us argue that the total space $F$ of the fibration in
\eqref{eq:fib-seq-with-n-1} is 
$n$-excisive, since the base and the fiber are. The argument is
straightforward but does require a variant of the five-lemma on $\pi_0$ to
make a conclusion about $\pi_0 F$. 
Once we know $F$ is $n$-excisive, we know $F \simeq P_n F$, as desired.
\end{proof}

\begin{corollary}
\label{cor:perp-F-zero}
If $F$ is a reduced functor that has either connected values
(\ref{hypothesis-1}) or 
group values (\ref{hypothesis-2}) on coproducts of $X$, and
$\Perp_{n+1} F(X) \simeq 0$, then $F(X) \simeq \Pnalg F(X)$.
\end{corollary}
\begin{proof}
$\Pnalg F(X)$ is defined by evaluating the functor $P_n (L F_X)$
at the space $S^0$, and $LF_X$ commutes with realizations (by
Lemma~\ref{lem:LF-commutes-with-realization}), so by
Proposition~\ref{prop:Pnalg-agrees-Pn}, we have
$P_n(L F_X)(S^0) \simeq P_n F(X)$. Hence 
the result follows from
Proposition~\ref{prop:perp-F-zero-realizations}.
\end{proof}


\section{Main Theorem, Case II: General $F$}
\label{sec:perp-F-nonzero}

In this section, the goal is to establish the other side of the
``ladder induction'' for Theorem~\ref{thm:main-theorem}.
We proceed essentially as outlined in
Section~\ref{sec:main-theorem-outline}. We actually decompose the
problem further, considering functors to discrete groups or connected
spaces. As one might expect, the case of a functor to discrete groups
is the pivotal one.


\subsection{Functors To Discrete Groups: $\Perp G^{ab} = 0$}

This section shows that a certain functor $G_n^{ab}$ has no $(n+1)^{\text{st}}$
cross effect, as one would expect in view of the construction
of $G_n^{ab}$ (given below).
In this section, we consider a functor $G$ from spaces to discrete
groups (for example, $G(X) = \pi_0 \Omega(X)$). 



\begin{definition}
\label{def:Gprime-Gab}
Given an $n>0$ and a functor $G$ to discrete groups,
define $G'_n := \Image (\epsilon: \Perp_{n+1} G
\rightarrow G)$ and $G^{ab}_n := \coker (\epsilon)$.  Usually, the $n$ is
clear from context, and we will abbreviate these $G'$ and $G^{ab}$.
\end{definition}
There is a short exact (fibration) sequence of groups
\begin{equation}
\label{eq:ses-Gprime-G-Gab}
G'(X) \rightarrow G(X) \rightarrow G^{ab}(X) ,
\end{equation}
and this sequence is surjective on $\pi_0$ (\emph{i.e.}, right exact). 

\begin{remark}
Note that since $\Perp_{n+1} G$ is constructed by taking a kernel,
the image of $G'$ is normal in $G$.
\end{remark}

Our motivation for the preceding notation comes from considering 
the case
when the source and target category under consideration are
both the category of groups and the functor $G$ is the identity $G(H)=H$.
In this case, the image of $\Perp_2 G(H)$ in $G(H)$ 
is the first derived subgroup of $H$.
The cokernel of the map $\Perp_2 G(H) \rightarrow G(H)$ is the
abelianization, $H^{ab}$. Lacking a more appropriate name for modding
out by higher derived subgroups, we continue to use the same notation
in that case.

\begin{lemma}
\label{lem:perp-Gab-0}
  If $F$ takes values in discrete groups, then  
  with $G'$ and $G^{ab}$ as in Definition~\ref{def:Gprime-Gab},
  $\Perp G^{ab}(X) \simeq 0$.
\end{lemma}
\begin{proof}
From the construction of $G'$, the map $\epsilon: \Perp G \rightarrow
G$ factors through the inclusion $i: G'\rightarrow G$. 
Applying $\Perp$ again results in the following diagram:
$$\xymatrix{
\Perp^2 G(X)
\ar[dr]_{\Perp \epsilon}
\ar[d]
&
\\
\Perp G'(X)
\ar[r]^{\Perp i}
&
\Perp G(X) 
\ar@/_1pc/[ul]_{\delta}
}$$
The map $\Perp \epsilon$ has a section, $\delta$, so it is surjective.
Hence $\Perp i$ is also surjective.
Consider the short exact sequence of functors to discrete groups
in Equation \eqref{eq:ses-Gprime-G-Gab} defining $G^{ab}$.
If we show that $\Perp$ preserves this short exact sequence, then the
surjectivity of the map $\Perp i$ will imply that $\Perp G^{ab} = 0$.

Short exact sequences of discrete groups are fiber sequences that are
surjective on the base space.
The functor $\Perp$ preserves fiber sequences because the construction
involves only taking fibers.
The functor $\Perp$ preserves surjections because all of the maps in
the cube $\CR_{n+1} F(X,\ldots,X)$ defining $\Perp_{n+1} F(X)$ 
have sections, and hence taking fibers with respect to them does not
lower connectivity.
\end{proof}

\subsection{$\Pnalg$ preserves $A_F$ fibration}

This section establishes that $\Pnalg$ actually produces a fibration
when applied to the fibration defining $A_F$.
The results in this section also contain a statement about the map
from $F \rightarrow \Pnalg F$, because in the case of $F$ taking
values in discrete groups, the proof that this map is surjective on
$\pi_0$ uses the same technical details that the proof that we get a 
fibration.

To remind the reader that the functor takes values in
discrete groups in the next proposition, we use the letter $G$ (for
group) to denote the functor, instead of the usual $F$.

\begin{proposition}
\label{prop:perp-f-nonzero-fib-discrete}
  If $G$ takes values in discrete groups (so in
  particular $G$ satisfies Hypothesis~\ref{hypothesis-2}), 
  then the following is a fibration sequence up to homotopy:
\begin{equation}
\label{eq:Pnd-still-fibration}
\Pnalg A_G(X) 
\rightarrow 
\Pnalg \left( \realization{\Perp_{n+1}^{*+1} G(X)}  \right)
\rightarrow 
\Pnalg G(X)
.
\end{equation}
\end{proposition}
\begin{proof}
Replacing the base $G$ in the definition of $A_G$
(Equation~\eqref{eq:AF-def}) with $G'$ from 
Definition~\ref{def:Gprime-Gab}, 
we have the fibration sequence
\begin{equation}
\label{eq:fib-Gprime-base}
A_G(X) 
\rightarrow 
\realization{\Perp_{n+1}^{*+1} G(X)}
\rightarrow
G'(X),
\end{equation}
and this sequence is surjective on $\pi_0$. 
The hypotheses of Corollary~\ref{cor:Pnd-group-connected-base} 
are satisfied by the
sequences in \eqref{eq:ses-Gprime-G-Gab} and \eqref{eq:fib-Gprime-base},
so applying $\Pnalg$ both are fibration sequences whose maps to the
base spaces are surjective on $\pi_0$:
\begin{gather}
\label{eq:Pnd-Gprime-G-Gab}
\Pnalg G'(X) 
\rightarrow 
\Pnalg G (X) 
\rightarrow 
\Pnalg G^{ab}(X)
\\
\label{eq:Pnd-AG-seq}
\Pnalg (A_G)(X) 
\rightarrow 
\Pnalg (\realization{\Perp_{n+1}^{*+1} G(-)}) (X)
\rightarrow
\Pnalg G' (X).
\end{gather}

The aim now is to show that \eqref{eq:Pnd-AG-seq} remains a fibration
when the base $\Pnalg G'(X)$ is replaced by $\Pnalg G(X)$.
From Lemma~\ref{lem:perp-Gab-0}, $\Perp_{n+1} G^{ab}(X) \simeq 0$, so
Corollary~\ref{cor:perp-F-zero} shows that $\Pnalg G^{ab}(X) \simeq 
G^{ab}(X)$, which is a discrete space. Then, using the long exact
sequence on homotopy, the fibration in \eqref{eq:Pnd-Gprime-G-Gab} gives
$\Pnalg G'(X) \xrightarrow{\simeq} \Pnalg G(X)$ except possibly on
$\pi_0$, where the map is injective. This is
enough to show that  changing the base in
\eqref{eq:Pnd-AG-seq} from $\Pnalg G'(X)$ to $\Pnalg G(X)$ still yields
a fibration. 
That is, \eqref{eq:Pnd-still-fibration} is a
fibration (but perhaps not surjective on $\pi_0$).
\end{proof}

\begin{proposition}
\label{prop:perp-f-nonzero-pi0-control}
If $G$ takes values in discrete groups (so in
particular $G$ satisfies Hypothesis~\ref{hypothesis-2}), then 
$$
\pi_0 
\Pnalg G(X) 
\cong \coker_{Gps} 
\left( \pi_0 \epsilon \right)
,$$
and the map $\pi_0 G \rightarrow \pi_0 \Pnalg G$ is the universal map
to the cokernel of $\pi_0 \epsilon$ in the category of groups.
\end{proposition}
\begin{proof}
As in the preceding
Proposition~\ref{prop:perp-f-nonzero-fib-discrete}, we have the following
fibration sequence that is surjective on $\pi_0$:
$$\Pnalg (A_G)(X) 
\rightarrow 
\Pnalg (\realization{\Perp_{n+1}^{*+1} G(-)}) (X)
\rightarrow
\Pnalg G' (X).
$$
Lemma~\ref{lem:Pnd-perp-contractible} shows that the total space
in this fibration is contractible, and the map to
the base is surjective on $\pi_0$, so  $\pi_0 \Pnalg G'(X) = 0$.

Also following Proposition~\ref{prop:perp-f-nonzero-fib-discrete}, 
we have the
following diagram in which the horizonal rows are fibrations that are
surjective on $\pi_0$:
$$\xymatrix{
G'(X)
\ar[r]\ar[d]
&
G(X)
\ar[r]\ar[d]
&
G^{ab}(X)
\ar[d]^{\simeq}
\\
\Pnalg G'(X)
\ar[r]
&
\Pnalg G(X)
\ar[r]
&
\Pnalg G^{ab}(X)
}$$
Since $\pi_0 \Pnalg G'(X) = 0$, the long exact sequence for the bottom
fibration implies that $\pi_0 \Pnalg G(X) \cong \pi_0 \Pnalg G^{ab}$.
The right hand vertical map is an equivalence, again as noted in the
preceding proposition, 
using  Lemma~\ref{lem:perp-Gab-0}  and Corollary~\ref{cor:perp-F-zero}.
Combining these, we have
\begin{align*}
\pi_0 \Pnalg G(X) 
&\cong 
\pi_0 \Pnalg G^{ab}(X)
\\
&\cong
\pi_0 G^{ab}(X)
,
\end{align*}
which is isomorphic to $\coker_{Gps} (\pi_0 \epsilon)$ because 
the map $\epsilon: \Perp_{n+1} G(X) \rightarrow
G(X)$ factors through $G'(X)$.
\end{proof}

\begin{proposition}
\label{prop:perp-f-nonzero-fib-conn}
  If $F$ has connected values (\ref{hypothesis-1}), 
 then the following is a fibration sequence up to homotopy:
\begin{equation*}
\Pnalg A_F(X)
\rightarrow 
\Pnalg \left( \realization{\Perp_{n+1}^{*+1} F(X)}  \right)
\rightarrow 
\Pnalg F(X)
,
\end{equation*}
and furthermore the map  $F(X)\rightarrow \Pnalg F(X)$ is (trivially)
surjective on $\pi_0$. 
\end{proposition}
\begin{proof}
If $F$ has connected values on coproducts of $X$, then
$$ A_F(X) 
\rightarrow 
\realization{\Perp_{n+1}^{*+1} F(X)}
\rightarrow 
F(X)
$$ 
is a fibration over a connected base. Therefore, by 
Corollary~\ref{cor:Pnd-group-connected-base}, applying $\Pnalg$ yields a
fibration, so 
Equation~\eqref{eq:Pnd-still-fibration} is a fibration.

The map $0\rightarrow F$ is $0$-connected, so
Theorem~\ref{thm:Pnd-preserves-connectivity} 
shows that $0 \simeq \Pnalg(0) \rightarrow \Pnalg F(X)$ is $0$-connected
as well. Hence $\pi_0 \Pnalg F(X) = 0$.
\end{proof}

\subsection{If $m<n$, then $\Pnalg \sR{\psp[n]{F}} \simeq 0$}

This section establishes the relatively easy fact that for $\Perp_n$, 
the part of the Goodwillie tower below degree $n$ is trivial.

\begin{lemma}
\label{lem:Pnd-perp-contractible}
  Let $R(X_1,\ldots,X_n) =  \realization{cr_n \left( \Perp_{n}^{*}
      F\right) (X_1, \ldots, X_n)}$ be a functor of $n$ variables. 
  Define the diagonal of such a functor to be the functor of one
  variable given by
  $(\diag R)(X) =
  R(X,\ldots,X)$.
  Then $P_m^{\alg} \left( \diag R \right) (X) \simeq 0$ for $0\le m<n$.
\end{lemma}
\begin{proof}
  Two results in Goodwillie's work \cite[Lemmas~3.1 and~3.2]{Cal3}
  combine to show that if 
  $H(X_1,\ldots, X_n)$ is a functor of $n$ variables that is
  contractible whenever some $X_i$ is contractible (this is called a
  ``multi-reduced'' functor), then $P_m (\diag
  H) \simeq 0$ for $0\le m<n$. 
  Writing out $P_m^{\alg} (\diag R) (X) = P_m [ L (\diag R)_X ] (S^0)$, 
  it is easy to check that
  $L(\diag R)_X$ is the diagonal of a multi-reduced functor, so
  Goodwillie's result applies. 
\end{proof}

\subsection{The functor $A_F$ has no $n+1$ cross effect}
Having created the functor $A_F$ to be ``$F$ with the 
cross effect killed'', we now need to establish that $\Perp A_F \simeq
0$. The main issue is the commuting of the $\Perp$ and the
realization. 

\begin{lemma}
\label{lem:perp-AF-zero}
Let $F$ be a functor that has either connected values
(\ref{hypothesis-1}) or group values
(\ref{hypothesis-2}), let $\Perp$
denote $\Perp_n$ for some $n$, and let $A_F$ be the functor given in
Definition~\ref{def:AF}.
Then $\Perp A_F$ is contractible.
\end{lemma}
\begin{proof}
Taking cross effects is a homotopy inverse limit construction, and
homotopy inverse limits commute, so
$$
\Perp A_F
\simeq 
\hofib \left(
\Perp \realization{\Perp^{*+1} F}
\rightarrow 
\Perp F
\right) .
$$
It is easy to check that if $F$ has connected values, then so does
$\Perp F$, and also $\realization{\Perp^{*+1} F}$.
If $F$ has group values, then so does $\Perp F$; hence $\Perp^{*+1} F$
is a functor to $\text{Simp}(\cat{G})$ --- simplicial grouplike $H$-spaces.
As remarked prior to Hypothesis~\ref{hyp:group-values}, the rigidity
of our definition of $H$-space implies the realization of
$\Perp^{*+1} F$ is still a functor to $\cat{G}$; that is,
$\realization{\Perp^{*+1} F}$ has group values.
By
Lemma~\ref{lem:perp-commutes-with-realization-for-simplicial-functors},
$\Perp$
commutes with the realization in that functor, so 
$\Perp\sR{\psp[]{F}} \simeq \sR{\Perp \psp[]{F}}$. Finally, the
existence of $\delta: \Perp F \rightarrow \Perp^2 F$ 
shows that $\Perp F$ is the augmentation of the simplicial space
$\Perp \psp[]{F}$, so the standard ``extra degeneracy'' argument
\cite[Exercise~8.4.6, p.~275]{Weibel:homological-algebra} shows that
$\sR{\Perp \psp[]{F}} \simeq \Perp F$, and hence that 
$$A_F \simeq
\hofib \left(
{\Perp F}
\rightarrow 
\Perp F
\right) \simeq 0, $$
as desired.
\end{proof}

\subsection{Proof of Main Theorem, Case II}
\label{sec:proof-of-perp-nonzero}

\begin{proposition}
\label{prop:perp-F-nonzero}
If $F$ is a reduced functor that
 has either connected values (\ref{hypothesis-1}) or
group values (\ref{hypothesis-2}) on coproducts of $X$, 
then the following is a fibration sequence up to homotopy:
\begin{equation}
\label{eq:perp-fibration-perp-nonzero-v2} 
\realization{\psp{F}(X)}
\xrightarrow{\epsilon}
F(X)
\rightarrow
\Pnalg F(X)
.
\end{equation}
Furthermore, the map 
$\pi_0 F(X) \rightarrow \pi_0 \Pnalg F(X)$
is the universal map to the cokernel of the group homomorphism 
$\pi_0 \Perp_{n+1} F(X) \rightarrow \pi_0 F(X)$.
\end{proposition}

\begin{proof}
First, suppose that $F(X)$ takes either connected values or discrete group
values on coproducts of $X$. 
Consider the auxiliary diagram created by applying $\Pnalg$ to the
fibration sequence defining $A_F(X)$:
$$\xymatrix{
A_F(X)
\ar[r]
\ar[d]
&
\realization{\Perp_{n+1}^{*+1} F(X)}
\ar[r]^{\epsilon}
\ar[d]
&
F(X)
\ar[d]
\\
\Pnalg A_F(X)
\ar[r]
&
\Pnalg \left( \realization{\Perp_{n+1}^{*+1} F(X)} \right)
\ar[r]
&
\Pnalg F(X)
}$$
Proposition~\ref{prop:perp-f-nonzero-fib-conn} (in the case of
connected values) or 
Proposition~\ref{prop:perp-f-nonzero-fib-discrete} (in the case of
discrete group values)
shows that the bottom row is
a fibration sequence up to homotopy. 
Proposition~\ref{prop:perp-f-nonzero-fib-conn} (connected values)
or 
Proposition~\ref{prop:perp-f-nonzero-pi0-control} (discrete group values)
imply that the map $F(X) \rightarrow \Pnalg F(X)$ surjective
on $\pi_0$. Lemma~\ref{lem:perp-AF-zero} shows that 
$\Perp_{n+1} A_F(X) \simeq 0$, so that
Corollary~\ref{cor:perp-F-zero} gives $A_F(X) \simeq \Pnalg
A_F(X)$, and hence the square on the right is homotopy Cartesian. 
Lemma~\ref{lem:Pnd-perp-contractible} shows that 
$\Pnalg \left( \realization{\Perp_{n+1}^{*+1} F(X)}\right) \simeq 0$,
so this square being Cartesian is equivalent to
\eqref{eq:perp-fibration-perp-nonzero-v2} 
being a fibration sequence up to homotopy, as we wanted to establish.

We can reduce the general problem when $F$ has group values
to the cases of connected and
discrete group values that we have already considered by examining
the fibration
$$ \widehat{F}(X) \rightarrow F (X) \rightarrow \pi_0 F(X),$$
where $\widehat{F}(X)$ is the component of the basepoint in $F(X)$.
This gives rise to the following square:
$$\xymatrix{
\realization{\Perp_{n+1}^{*+1} \widehat{F}(X)}
\ar[r]
\ar[d]
&
\widehat{F}
\ar[r]
\ar[d]
&
\Pnalg \widehat{F}
\ar[d]
\\
\realization{\Perp_{n+1}^{*+1} {F}(X)}
\ar[r]
\ar[d]
&
{F}
\ar[r]
\ar[d]
&
\Pnalg {F}
\ar[d]
\\
\realization{\Perp_{n+1}^{*+1} \pi_0{F}(X)}
\ar[r]
&
\pi_0{F}
\ar[r]
&
\Pnalg \pi_0 {F}
}$$
It is straightforward to check that every row and column except the
middle row is a fibration that is surjective in $\pi_0$, and that the
composition of the two maps in the middle is null homotopic, and 
we have shown that the map $F \rightarrow \Pnalg F$ is surjective on
$\pi_0$. 
This gives us the data required to use the $3 \times 3$ lemma for
fibrations to show that the 
middle row (\emph{i.e.}, \eqref{eq:Pnd-still-fibration}) is a
fibration and surjective on $\pi_0$.

The statement about $\pi_0$ is trivial in the connected case;
$\pi_0$ of every space in the top row is zero (which is trivially a
group). This implies that the vertical arrows connecting the second
and third rows are $\pi_0$-isomorphisms, so the statement about
$\pi_0$ follows from Proposition~\ref{prop:perp-f-nonzero-pi0-control}.
\end{proof}



\bibliographystyle{amsplain}
\bibliography{andrew-math}

\end{document}